\documentclass[12pt, twoside]{report}

\usepackage{amsfonts}
\usepackage{amssymb}
\usepackage{amsmath, amsthm, latexsym}

\def \hf{\hat{F}}
\def \c{\mathbb{C}}
\def \z{\mathbb{Z}}
\def \r{\mathbb{R}}
\def \p{\mathbb{P}}
\def \o{\mathbb{O}}
\def \j{\mathbb{J}}

\def \tor{{\c^{*}}^n}
\def \Tor{{\c^{*}}^N}
\def \t{\mathfrak{t}}
\def \k{\mathfrak{k}}
\def \g{\mathfrak{g}}
\def \b{\mathfrak{b}}
\def \SL{SL_2(\c)}
\def \.{\cdot}

\theoremstyle{plain}
\newtheorem{Th}{Theorem}[section]

\newtheorem{T}{Theorem} 

\newtheorem{Lem}[Th]{Lemma}
\newtheorem{Prop}[Th]{Proposition}
\newtheorem{Cor}[Th]{Corollary}

\theoremstyle{definition}
\newtheorem{Ex}{Example}[section]
\newtheorem{Def}{Definition}[section]
\newtheorem{Rem}{Remark}[section]
\newtheorem{Not}{Notation}[section]

\begin{document}

\title{Morse theory and Euler characteristic of sections of 
spherical varieties}
\author{Kiumars Kaveh\\
\\
Department of Mathematics\\
\\
University of Toronto}
\maketitle

\tableofcontents

\chapter{Introduction}

A classical result by D. Bernstein (cf. ~{Bernstien}) 
asserts that if $p(x_1, \ldots,
x_n) = \sum_{\alpha \in \mathbb{Z} ^n}
c_{\alpha} x^{\alpha}$, (where $\alpha = (\alpha_1,
\ldots, 
\alpha_n)$ and
$x^\alpha = (x_{1}^{\alpha_1}, \ldots, 
x_{n}^{\alpha_n}))$ is a Laurant polynomial
with generic coefficients then ,
Euler characteristic (in homotopy sense) of the
hypersurface $\{ x in \c^n| p(x) = 0\}$ is
equal to $n! Vol(\Delta)$
where $\Delta$ is the Netwon polyhedron of $p$ i.e. the convex
hull of points $\alpha \in \mathbb{Z}^n$ with $c_\alpha
\neq 0$. 
One can put this result in a more fancy way as
follows: 

Let $p(x) = \sum_{\alpha \in S}c_\alpha x^\alpha$ with 
$S = \{ \alpha_1, \ldots, \alpha_N \} \subset \mathbb{Z}
^N$ and let $\pi : \tor \rightarrow \Tor \subset GL(N,
\mathbb{C})$
be
the finite dimensional representation of algebraic torus
$\mathbb{C}^{*^n}$ defined by $\pi(x) = diag(
x^{\alpha_1}, \ldots, x^{\alpha_N})$ and consider
the hyperplane  $L = \{ x \in \c^n | f(x) = 0 \}$ where $f(x_1,
\ldots, x_N) = \sum_{i=1}^{N} \alpha_i x_i$, then if
$\alpha_i$ are generic enough, the Euler characteristic of 
hyperplane section $L \cap \pi(\mathbb{C}^{*^n})$
 is equal
to
$n! Vol(\Delta)$, where $\Delta$ is convex hull of
$\{\alpha_1, \ldots, \alpha_N \}$.

In this thesis, we try to generalize the above result to actions of 
other algebraic groups. 

On the other hand, there is 
the famous theorem of Kuchnierenko which states that the degree of $\pi(\tor)$
sitting inside $\Tor$ is equal to $n! Vol(\Delta)$. Hence Bernstein's
result indeed claims that upto a sign, Euler characteristic of a generic
hyperplane section of $\pi(\tor)$ is equal to $deg(\pi(\tor))$ as a subvariety
of $\c^N$. 
 
We will use some variant of Morse theory to prove our results. For the 
manifold we take a submanifold of $\c^N$ and for the Morse function we 
take a linear functional on $\c^N$. This is different from the classical 
Morse theory since neither the manifold is assumed to be compact nor the 
function is proper. 

Morse theory relates the Euler characteristic of sections and the number of 
critical points. Let us explain the idea: 

First consider the real case. Let $f$ be a Morse function from a compact manifold $M$ 
to $\r$. From Morse theory we know how the topology of the sets $M_{\leq a}= f^{-1}((-\infty, a])$ 
changes as $a$ passes a critical value. More precisely, suppose $c$ is the only critical 
value between $a$ and $b$, $a \leq b$ and $f^{-1}(c)$ contains only one critical point $p$. 
Then $M_{\leq b}$ has the homotopy type of $M_{\leq a}$ with a cell
of some dimension $\lambda$ attached. $\lambda$, called the \textit{index of the critical point $p$}
, is in fact equal to the number of negative eigen values of the matrix of 
second derivative of $f$.

As for the complex case, let $f: M \rightarrow \c$ be a Morse function from a 
$d$-dimensional complex manifold $M$ 
to $\c$. We do not assume that $M$ is compact (as there are no non-constant holomorphic functions
on compact manifolds). But suppose that we can apply Morse theory to the real part 
of $f$.
If $f$ does not have any critical points, Morse theory tells us that $f$ defines a fibration of 
$M$ over $\c$, that is
for any value of $c$ in the range of $f$ we have
$$ M \cong  f^{-1}(c) \times \c .$$
and hence for the Euler characteristic $\chi$ (in the homotopy sense) we have: 
$$\chi(M) = \chi(f^{-1}(c)) \cdot \chi(\c) = \chi(f^{-1}(c))$$ , as $\c$ has the 
same homotopy type as a point and $\chi$ of a point is $1$. 
But usually functions do have critical points. In this case, away from the 
set of critical points $f$ is a fibration. Moreover Mores theory tells us
about the topology of fibres at critical values: roughly speaking,
as one moves from a 
regular value to a critical value, the topology of the fibre changes in such a way 
that a sphere of real dimension $d$ vanishes to a point. So in terms of 
Euler characteristic we have:
$$ \chi(M) = \chi(f^{-1}(c)) + \text{correction terms coming from critical 
points}.$$
where $c$ is a regular value. 
As the Euler characteristic of a punctured sphere of real dimension 
$d$ is $(-1)^{d}$, we get $$ \chi(M) = \chi(f^{-1}(c)) + (-1)^{d} 
\cdot \text{the number of critical points}.$$

Let us denote the number of critical points of $f$ on $M$ by $\mu(M, f)$.
We can formulate the above result as 
\begin{Th}[Formula for Euler characterisitic]
Let $M$ be a (closed) complex algebraic submanifold of $\c^n$ 
of complex dimension $d$, and $f$ a generic complex linear functional 
on $\c^n$. Let $c$ be a regular value for $f_{|M}$, we have
$$ \chi(f^{-1}(c)) = \chi(M) + (-1)^{d+1} \cdot \mu(f, M).$$   
\end{Th}
 
The thesis is made up of two parts. The first part deals with Morse theory. It
developes a variant of Morse theory for $f: M \rightarrow \r$ where $f$ is not 
necessarily proper and $M$ is an algebraic submanifold of $\r^N$ such that 
$f$ satisfies certain transversality condition with respect to $M$. This 
can be roughly interpreted as $f$ does not have any critical point at infinty of $M$.

From this Morse theoretic results we can then derive a formula for the Euler characteristic of 
hyperplane sections of algebraic submanifolds of $\c^N$, in terms of number of critical 
points of the linear functional defining the hyperplane.

We will approach this variant of Morse theory in several ways:

1. Stratified Morse theory: We will derive our required result from the 
Thom-Mather-Whitney stratification theory. They basically generalized the 
first theorem of calssical Morse theory to the stratified spaces. \footnote{
The generalization of the second theorem of Morse theory to stratified space was 
done by Goresky and McPherson, cf. \textit{Stratified Morse Theory} {}.}

2. Generalized Morse theory of Palais and Smale: One can also derive our
result from the Generalized Morse theory of Palais and Smale. Their 
theory is a generaliztion of classical Morse theory to the Hilbert manifolds.
Main application of their theory is for infinite dimensional manifolds of geodesics
(cf. \textit{Generalized Morse theory} {} or \textit{Morse thoery on Hilbert manifolds}
{}). Interesting enough this was also the original motivation of Morse himself.
  
3. Finally, in case our manifold is an orbit of a Lie group action we will 
give direct proofs for the results. In fact as will be explained in its place,
our results obtained are stronger than the results for a general stratification:
\textit{for orbits of Lie group actions we can prove that the 
sections are diffeomorphic (rather than only homeomorphic).}

The main objective of the second part is to apply the formula for 
the Euler characteristic obtained in the previous part to prove a generaliztion 
of Bernstein's theorem to the so-called \textit{actions with spherical orbits}.

We start by discussing the classical results for torus actions, Kouchnierenko's theorem 
and its generalization to spherical varieties.

Next we attempt to generalize Bernstein's Euler characteristic theorem to representations
of reductive algebraic groups (Section ~{}). Let $G$ be a complex connected reductive 
algebraic group and $\pi: G \rightarrow GL(n,\c)$ a faithful representation.
To apply Morse theory we need $\pi(G)$ be a closed in $M(n, \c)$, the vector space of 
$n \times n$ matrices. We prove a criterion for $\pi(G)$ to be closed (Proposition ~{})
and from the formula for Euler characteristic we then obtain

\begin{Th}
Let $\pi: G \rightarrow GL(N, \c)$ be a faithful representation of
a $d$-dimensional complex connected reductive group $G$. Suppose
$\pi(G)$ is closed in $M(N, \c)$ that is origin belongs to the
convex hull of weights of $\pi$. Then for $f$ a generic linear
functional on $M(N, \c)$ and $c$ a generic complex number we have
$$ \chi(\{x \in \pi(G)| f(x) = c \}) = (-1)^{(d+1)} \cdot \mu(f,
\pi(G)).$$
\end{Th}
   
We can then prove Bernstein's Euler characteristic theorem by showing that the degree 
(as a subvariety) of a generic orbit of $\tor$ acting on $\c^N$ is equal to the number 
of critical points of a generic functional restricted to this orbit. \footnote{Bernstein 
himself, proved his result using similar methods, cf. {Bernstein}}
Next we investigate in what other group actions one has that the number of critical points on
an orbit is equal to the degree. We will carefully consider the case of representations of 
$\SL$ and see that the number of critical points on an orbit is NOT equal to degree in this case. 
In fact, it turns out that equality of number of critical points and degree is a special 
situation and is usually not true. Our main result will be that it is indeed true for linear 
actions with so called \textit{generic spherical orbits}. These kind of actions already contain 
all linear torus actions. 

\begin{Th}
Let $G$ acts linearly on $\c^N$ such that generic orbits are
spherical and closed. Let $X$ be a generic orbit and let $f$ be a
generic linear functional on $\c^N$. Then
$$ deg(X) = \mu(X, f).$$
\end{Th}
\begin{Cor}[the main theorem]
Let $G$ acts linearly on $\c^N$ such that generic orbits are
spherical and closed. Let $X$ be a generic orbit and let $f$ be a
generic linear functional on $\c^N$ and $c$ a generic complex
number. Then
$$ \chi(f^{-1}(c) \cap X) = \chi(X) + (-1)^{dim(X)+1} \cdot
deg(X).$$
\end{Cor}
I should mention that, being a homogeneous space of $G$, $\chi(X)$ is usually zero and 
in case it is not zero there is a simple formula for $\chi(X)$ in terms of stabilizer 
subgroup of $X$ and the Weyl group of $G$ (see Proposition ~{}).
  
In section ~{} we give the proof of the above result. 
Interesting enough, all representations with generic spherical orbits had already been
classified by I. Arzhantsev. His list of \textit{indecomposable actions with spherical orbits}, 
includes all torus actions as well as about 30 
more examples. We will examine all the examples and verify again directly 
that the degree is equal to the number of critical points.
Main concrete example will be:

\begin{Th}
The Euler characteristic of a generic hyperplane section of 
$SL(n, \c) \subset M(n, \c)$ is equal to ${-1}^{n} \. n$.
\end{Th}

Finally, in the last section, we give a formula, in general, for the number of critical points 
in terms of degree and the intersection numbers of Chern classes. The number of critical 
points would be then equal to degree if all the terms corresponding to intersection numbers 
Chern classes cancel each other out.  

\chapter{Morse theory}

\section{Basic definitions and the classical Morse theory}
In this section we briefly go over the basic definitions 
and state the two classical theorems of Morse theory.

Through out this section, $M$ is a $d$-dimensional smooth manifold and 
$f: M \rightarrow \r$ is a smooth function. 

\begin{Def}
A $p \in M$ is called a \textit{critical point} of $f$ if $df(p) = 0$. 
More generally if $f: M \rightarrow N$ is a differentiable map between smooth manifolds,
a point $p \in M$ is called a critical point of $f$
if the derivative $df(p): T_P M \rightarrow T_{f(P)} N$  is not surjective.   
Points which are not critical are called \textit{regular}. 
\end{Def}

\begin{Def}
A $c \in \r$ is called a \textit{critical value} if $f^{-1}(c)$ contains a 
critical point, otherwise $c$ is called a \textit{regular value}.
A critical value which has only one inverse image is 
called a \textit{simple critical value}. 
\end{Def}

The celeberated theorem of Sard {Milnor}, says that almost every
value is regular, i.e. the set of critical values is of measure zero.

The second derivative $d^2 f$ at each point $p \in M$ is a bilinear 
map on $T_pM \times T_pM$. If one fixes a local coordinate system 
at $p$ on M, this bilinear map is represented by the $d \times d$ matrix 
$[ \partial f / \partial x_i \partial x_j (p)]$ of second order 
partial derivatives of $f$. This matrix is calles the \textit{Hessian 
matrix of $f$ at $p$} with respect to the local coordinate system chosen.

\begin{Def}
A critical point $p \in M$ is called \textit{non-degenerate}, if $d^2 f$
is a non-degenerate bilinear form, or equvalently if the Hessian 
matrix, with respect to some coordinate system,is invertible. A function 
, all whose critical points are non-degenerate is called a \textit{Morse
function}. If all the critical values of $f$ are simple $f$ is called a 
\textit{simple Morse function}.
\end{Def}

\begin{Def}
Let $p \in M$ be a non-degenerate critical point. The number of 
negative eigen values of the Hessian matrix of $f$ at $p$, in some 
local coordinate, is called the \textit{index of the critical point $p$}. 
One can easily see that this number is independent of the local coordinate 
chosen.
\end{Def}

It can be proved that the set of simple Morse functions is a dense open 
set in the space of all $C^{\infty}$ functions. Any function after an
arbitrarily small perturbation becomes a Morse function with simple 
critical values,

The key lemma in determination of local behaviour of functions at 
non-degenerate critical points is the so-called \textit{Morse lemma} 
(cf. ~{Milnor}). It is an essential step in the proof of Morse theorems.

\begin{Th}[Morse Lemma]
Let $f$ be smooth function on $\r^n$ such that origin $O$ is a non-degenerate 
critical point of $f$ and $f(O) = 0$.
Then with a smooth change of coordinates in a neighbourhood of the origin
$f$ will have the form 
$$ f(x_1, \ldots, x_n) = -{x_1}^2 - {x_2}^2 - \cdots - {x_k}^2 + 
{x_{k+1}}^2 + \cdots {x_n}^2.$$
where $k$ is the index of the critical point $O$.
\end{Th}
 
Now we are ready to state the Morse theorems. They are concerned with the
topology of the sets $M_{\leq a} = f^{-1}((-\infty, a])$, that is all the 
points whose $f$ is below $a$. The proof can be found for example in 
~{Milnor} or ~{Hirsch}.

\begin{T}[Classical Morse theory part A]
Let $M$ be a $C^{\infty}$ manifold and let $f: M \rightarrow \r$ be a 
Morse function. Suppose there are no critical values in the interval 
$[a, b]$. Then the subsets $M_{\leq a}$ and $M_{\leq b}$ have the same homotopy type.
In other words, as $c$ varies between the open interval between 
two adjacent critical values, the homotopy type of $M_{\leq}$ remains 
constant. 
\end{T}

\begin{T}[Classical Morse theory part B]
Let $M$ be a $C^{\infty}$ manifold and let $f: M \rightarrow \r$ be a 
simple Morse function. Let $c$ be the only critical value in the interval
$[a, b]$. Then $M_b$ has the homotopy type of $M_a$ with a cell of dimension
$\lambda$ attached, where $\lambda$ is the index of the critical point
$f^{-1}(c)$.   
\end{T}

We should remark that it is not crucial in the above theorem to assume
$c$ is a simple critical value. In general, we get one cell attached for each 
critical poin in the inverse image of $c$.

\section{Main theorems of the chapter}

\subsection{A variant of Morse theory}
Suppose $M \subset \r^n$ is a (closed) algebraic submanifold and 
$f: \r^n \rightarrow \r$ is a linear functional, such that $f_{|M}$ is a 
Morse function. In this section we formulate a Morse theory which 
says that for good enough $f$ (with respect to $M$) we have analogues of 
theorems A and B of the classical Morse theory for $f_{|M}$. \footnote{As we 
will see, $M$ being 
algebraic is not crucial. We only need that closure of $M$ in projective space 
admits a Whitney stratification with finite number of strata.}

The first part of Morse theory (theorem A) says that the homotopy type of the 
sets $M_{\leq c}$ remains constant as long as $c$ is varying between 
two adjacent critical points. Let us see an example where this is not true
if $M$ is not compact.

\begin{Ex}
Let $M \subset \r^2$ be the right part of the graph of the function $y = 1/x$.
and let $f(x, y) = y$. As you see $M_{\leq c} = \emptyset$ for $c \leq 0$ and 
$M_{\leq c}$ consists of a single point for $c > 0$, i.e. tha homotopy type of
$M_{\leq c}$ changes although $c$ is not a critical value for $f$. 
\end{Ex}

In the above example there is no critical point $p \in M$ with $f(p) = 0$
but the line $y = 0$ is tangent to $M$ at infinity suggesting that we can 
think of $c = 0$ be a critical value corresponding to a point at infinity of 
$M$.

The philosophy is that one can repeat the Morse theory for non-compact 
manifolds (or non-proper functions) if there are no critical points at
infinity. But of course, one should make it more precise what one means by
$f_{|M}$ has no critical point at infinity. 

We now state the theorems, later on we discuss the condition of 
having no critical point at infinity. As we will see 
almost all the linear functionals in $\r^n$ behave well and have no 
critical point at infinity for given algebraic submanifold $M$.

\begin{Th}[A$'$]
 Let $M$ be a (closed) algebraic submanifold of $\r^n$ and 
 $f$ a generic linear functional on $\r^n$. Suppose 
 $f_{|M}$ does not have any critical value in the interval $[a, b]$. 
 Then the set $M_{\leq a}$ and $M_{\leq b}$ have the 
 same homotopy type.
\end{Th}

From Theorem A$'$, one can prove a more general form of it,
i.e. when we have a function from $M$ to $\r ^k$. 
 
\begin{Th}[A$'$, the general form]
Let $M$ a (closed) algebraic submanifold of $\r^n$ and 
$f$ a generic linear function from $\r^n$ to $\r^k$.
Suppose $U \subset \r^k$ is an open set that contains no critical 
values. Then $f: f^{-1}(U) \cap M \rightarrow U$ is a fiberation.  
\end{Th}

As soon as theorem A$'$ is established, repeating the proof of
second part of the classical Morse theory (theorem B) we obtain

\begin{Th}[B$'$]
Let $M$ be a (closed) algebraic submanifold of $\r^n$ and 
$f$ a generic linear functional on $\r^n$
let $c$ be the only critical value in the interval
$[a, b]$. Then $M_b$ has the homotopy type of $M_a$ with a cell of 
dimension $\lambda$ attached, where $\lambda$ is the index of the critical 
point $f^{-1}(c)$.   
\end{Th}
\begin{proof}
Let $p$ be the critical point $f^{-1}(c)$.
Following the proof of Morse theorem (cf. Morse Theory, J. Milnor,
Theorem $3.2$, Chap.$1$, ~[Milnor]) one can find $\epsilon > 0$ 
and arbitrarily small neighbourhood $U$ of $p
\in \r^ {n}$ and $\tilde{f} : \r^ {n} 
\rightarrow \r$, such that:
$\tilde{f} = f$ outside $U$ and inside $U$,
$\tilde{f}$ is 
defined in the following way:

In $U$ consider the coordinate system $u^1, \ldots, u^{n}$ such that:
\begin{enumerate}
\item $u^1(p) = u^n(p) = 0$
\item $Y \cap U \subset \{ u^{m+1} = \cdots = u^n = 0 \}$
\item $f_{|_Y}(x) = c - (u^1)^2 - \cdots -(u^{\lambda})^2 +
(u^{\lambda
+1})^2 + \cdots + (u^m)^2$, where $m = dim(M)$. 

\end{enumerate}
Now define $\tilde{f}$ in $U$ by:
$$ \tilde{f} = f - \mu((u^1)^2 + \cdots + (u^\lambda)^2 + 2(u^{\lambda
+1})^2 + \cdots + 2(u^m)^2 + (u^{m+1})^2 + \cdots + (u^{n})^2)$$
where $\mu :\r \rightarrow \r$ is smooth and $\mu(0)>
\epsilon, \mu(r) = 0, \forall r\geq2\epsilon$ and $-1 < \mu'(r) \leq 0,
\forall r$. It is clear from defenition of $\tilde{f}$ that it is smooth
on whole $\r ^{2n}$.

As a function on $M$, $\tilde{f}$ has the following properties:
\begin{enumerate}
\item $\tilde{f}^{-1} (-\infty , c-\epsilon] = f^{-1} (-\infty,
c- \epsilon]$
\item critical points of $\tilde{f}$ and $f_1$ on $M$ are the same.
\item $\tilde{f}$ has no critical value in $[c-\epsilon, c+\epsilon]$.
\item (Main Property) $\tilde{f}^{-1} (-\infty, c+\epsilon]$ has 
homotopy type of $f ^{-1} (-\infty, c+\epsilon]$ with a cell 
of dimension $\lambda$ attached.
\end{enumerate}

Now $\tilde{f}$ satisfies the conditions of theorem A$'$ on 
$[c-\epsilon, c+\epsilon]$ as $f = \tilde{f}$ outside $U$. Hence we conclude 
that 
$$ \tilde{f}^{-1}(-\infty, c-\epsilon] \sim 
\tilde{f}^{-1}(-\infty, c+\epsilon]. $$
And from the property number 4 of $\tilde{f}$
\begin{eqnarray*}
f^{-1}(-\infty, c+\epsilon] 
&\sim& \tilde{f}^{-1}(-\infty, c+\epsilon]
\textup{ with a cell of dimension } \lambda \textup{attached}. \cr
&\sim& \tilde{f}^{-1}(-\infty, c-\epsilon]
\textup{ with a cell of dimension } \lambda \textup{ attached}. \cr
&=& f^{-1}(-\infty, c-\epsilon]
\textup{ with a cell of dimension } \lambda \textup{ attached}. \cr
\end{eqnarray*}
\end{proof}
 
\subsection{Formula for the Euler characteristic of sections}
In this subsection, we apply the above theorems 
to obtain a formula for the Euler characterisitic of hyperplane sections of 
algebraic submanifolds.

M will denote a (closed) complex algebraic submanifold of 
$\c^n \cong \r^{2n}$ of complex dimension $m$. To make $M$ compact, we
can either consider its closure in $\r P^{2n}$ or $\c P^n$. Since our Morse 
theoretic theorems A$'$ and B$'$ deal with real projective space, we 
prefer to take the closure of $M$ in $\r P^{2n}$ which we denote by
$\overline{M}$. By theorem ~{}, there is a stratification $\mathcal{A}$
of $\overline{M}$ with finite number of algebraic strata such that $M$ itself
is a union of strata.  

The definitions of critical point and 
non-degenerate critical point for the complex functions is verbatim to 
the real functions.

\begin{Not}
We denoe the number of critical points of a function $f$ from $M$ to $\r$ (or
$\c$)by $\mu(M, f)$ (whenever this is finite).
\end{Not}

We prove

\begin{Th}[Formula for Euler characterisitic]
Let $M$ be a (closed) complex algebraic submanifold of $\c^n$ 
of complex dimension $m$, and $f$ a generic complex linear functional 
on $\c^n$. Let $c$ be a regular value for $f_{|M}$, we have
$$ \chi(f^{-1}(c)) = \chi(M) + (-1)^{m+1} \cdot \mu(f, M).$$   
\end{Th}

\begin{proof}
Write $f$ as $f_1 + i \. f_2$. Obviously, $f_1$ and $f_2$ are $\r$-linear
functionals on $\c^n \cong \r^{2n}$. To prove the theorem we apply the 
theorems A$'$ and B$'$ to $f_1$. We need few observations about the 
critical points of $f$ and $f_1$. Note that by Cauchy-Riemann relations, 
$p$ is a critical point of $f_{|M}$ if and only if, it is a critical 
point of ${f_1}_{|M}$.

\begin{Th}[Complex Morse Lemma]
Let $f$ be a holomorphic function on $\c^m$ such that origin $O$ is a 
non-degenerate critical point of $f$ and $f(O) = 0$.
Then with a smooth change of coordinates in a neighbourhood of the origin
$f$ will have the form 
$$ f(x_1, \ldots, x_m) = {x_1}^2 + {x_2}^2 + \cdots + {x_m}^2.$$
\end{Th}
\begin{proof}
The proof is verbatim to the real Morse lemma. Only notice that 
any non-degenerate complex quadratic form, after a linear 
change of coordinates, can be put in the form ${x_1}^2 + \cdots + {x_m}^2$. 
\end{proof}

\begin{Cor}
Let $f$ be a holomorphic function on $U \subset \c^m$ and $f_1 = Re(f)$.
Then any non-degenerate critical point of $f_1$ (which is automatically 
a critical point of $f$) has index equal to $m$.
\end{Cor}
\begin{proof}
Simple calculation shows that the Hessian matrix for the function 
$x^2 - y^2$ on $\r^2$ (which is the real part of the function 
$z^2$) at origin is 
$\left[\begin{matrix}
2 & 0 \\
0 & -2 \\
\end{matrix}
\right] .$ i.e. there is one positive and one negative eigen value. 
From this we can easily see that the index of zero as a non-degenerate 
critical point of $Re({x_1}^2 + \cdots + {x_m}^2)$ is $m$.  
\end{proof}

Now since $f$ is an algebraic function on $M$
it has only finitely many critical points on $M$, also as we mentioned 
before, for generic $f$, all $\overline{\{ x \in \c^n | f(x)=c\}}$, 
$(\forall c \in \c)$, are transverse to all strata at infinity in
$\overline{M}$.

Take the regular values $c, d \in \r$ such that $f_1^{-1}(-\infty, c]$ 
does not contain any critical point  and $f_1^{-1}(-\infty
, d]$ contains all critical points of $f_1$ on $M$.
It follows easily from the above corollary of complex Morse lemma that 
the index of a critical point $p$ of the function $f_1$ on $M$ is 
equal to $m$, the complex dimension of the manifold.  

By repeated application of theorems A$'$ and B$'$ 
(for all the critical values) we get

\begin{eqnarray*}
M &=& f_1^{-1}(-\infty, +\infty) \cr
&\sim& f_1^{-1}(-\infty, d] \cr
&\sim& f_1^{-1}(-\infty, c] \text{ with cells of real dimension } d
\text{ attached }
\end{eqnarray*}
where $\sim$ means homotopy equivalent. 

In terms of the Euler characteristic 

\begin{eqnarray*}
\chi(M) &=& \chi(f_1^{-1}(-\infty,c]) + (-1)^m \. \mu(M, f_1). \cr 
&=& \chi(f_1^{-1}(-\infty,c]) + (-1)^m \. \mu(M, f). \cr
\end{eqnarray*}

But by Theorem A$'$ (general form), applied to $f$, $f_1^{-1}(-\infty, c]$ 
has the same homotopy type
as $f^{-1}(c)$ (because $f:M \rightarrow \c$ is a 
fiberation restricted to $f_1^{-1}(-\infty, c])$, so we get

$$\chi(M) = \chi(f^{-1}(c))+ (-1)^{m} \cdot \mu(f, M). $$
or 
$$\chi(f^{-1}(c)) = \chi(M)+ (-1)^{m+1} \cdot \mu(f, M). $$

\end{proof}

\section{Critical points at infinity, Whitney stratified sets
and the proof of theorem A$'$}

The genericity condition for $f$ in the above theorems is that 
$f$ has no \textit{critical points at infinity}. We are going to make this 
precise. A point being regular point can be stated as a transversality
condition: a point $p \in M$ with $f(p) = c$ is regular iff $H_c = 
\{ x \in \r^n | f(x) = c \}$ is transverse to $M$ at $p$. In a similar 
way one can say that $f_{|M}$ does not have any critical point at
infinity iff level sets of $f$ are transverse to $\overline{M}
\setminus M \subset 
\r P^n$. But $\overline{M}$ is not necessarily a manifold. in order to resolve
this problem we introduce notion of a stratified set and a Whitney 
stratified set.

\subsection{Whitney stratification}

In many places in math, we deal with objects such as algebraic varieties
which are not manifolds but are union of manifolds glued toether, stratification
theory provides a general frame work for doing analysis on this objects

\begin{Def}[Stratified set]
Let $Y$ be a closed subset of a smooth manifold $X$. A collection 
$\mathcal A$ of disjoint, locally closed submanifolds $ Si \subset Y (i\in I)$
is called a \textit{startification of $Y$} iff
\begin{enumerate}
\item $ Y = \cup_{i \in I} S_i$. 
\item $S_i \cap S_j \neq \emptyset$ only if $i=j$ or $S_i \subset \overline{S_j}
$ or $S_j \subset \overline{S_i}$.
\end{enumerate}

The submanifolds $S_i$ are then called \textit{strata} and $Y$ is called a 
\textit{stratified set}.
\end{Def}

Whitney realized that the notion of a stratification, in general, 
is too wild to expect a good theory for. One needs some 
conditions to ensure that the strata are glued together in a regular
way. To this end, he imposed his conditions (A) and (B) on a stratification 
to guarantee the so called \textit{topological trivality along the strata}
, i.e. if we slice the stratified set and move the slice along some 
strata, the topological picture does not change. 

\begin{Def}[Whitney stratified set]
Let $Y$ be a closed subset of $\r^n$. A stratification $\mathcal A$
of $Y$ is called a \textit{Whitney stratification} iff for any pair of 
strata $A$ and $B \in \mathcal A$ with $A \subset \overline{B}$ 
we have the following conditionds (A) and (B) satisfied

\begin{itemize}
\item[(A)] For any point $a \in A$ and a sequence $\{ b_i \} 
(i = 1, 2, \ldots )$ of points in B, assume that 
$\lim_{i \rightarrow \infty} b_i = a$ and there is a limit plane 
$L = \lim_{i \rightarrow \infty} T_{b_i}B \subset T_a \r^n$. Then 
we have $T_a A \subset L$.
\item[(B)] For any point $a \in A$ and sequences $\{ a_i \in A\}$ and 
$\{ b_i \in B\}$ assume that $\lim_{i \rightarrow \infty} a_i =
\lim_{i \rightarrow \infty} b_i = a$. Write $\overline{b_ia_i}$ for the 
secant line through $b_i$ and $a_i$. Think of $\overline{b_ia_i}$ as a 
subspace of $T_{b_i}\r^n$. Assume that there are limit plane 
$L = \lim_{i \rightarrow \infty} T_{b_i}B \subset T_a \r^n$ and 
limit line $K = \lim_{i \rightarrow \infty} \overline{b_ia_i} 
\subset T_a \r^n$. Then we have: $K \subset L$. 
\end{itemize}
\end{Def}

One can prove that Whitney conditions (A) and (B) are invariant 
under diffeomorphisms of $\r^n$. Hence we can speak of a Whitney
stratified subset of a smooth manifold.

It is easy to see that condition (B) implies (A) and hence, logically
speaking, it is redundant. But there are still reasons that authors prefer 
to mention both of them together.

Whitney stratification is important because one can prove 
(cf. ~{Goresky-McPherson} p.)

\begin{enumerate}
\item Whitney stratifications are locally topologically trivial
along the strata. 
\item Any closed analytic subset of of an analytic manifold 
admits a Whitney stratification.
\item Whitney stratified spaces can be triangulated.
\item \label{intersection-of-strat} The transversal intersection of two Whitney 
stratified spaces is again a Whitney stratified space, whose strata are the 
intersections of the strata of the two spaces. Also the Cartesian product of 
two Whitney stratified spaces is a Whitney startified space, the strata being the product 
of starta of the stratifications.
\end{enumerate}
and the most important example of Whitney stratification for us is
the case of algebraic varieties:
\begin{Th}[cf. ~{}]
Any (closed) algebraic subset $M$ of an algebraic manifold $X$ (either 
over $\c$ or $\r$) admits a Whitney stratification with finite number 
of algebraic strata. Moreover if we fix an algebraic subvariety $V$
of $M$, we can choose the stratification such that $V$ becomes a union of 
strata.
\end{Th}

Let $M \subset \r^n$ be a (closed) algebraic submanifold and let 
$\mathcal A = \{ X_0, X_1, \ldots, X_r \} (X_0 = M)$ be a finite 
Whitney stratification of $\overline{M} \subset \r P^n$. As usual 
let $f(x) = \sum_{i=1}^{n}f_ix_i$ be a linear functional on $\r^n$.
Consider $\r P^n$ as $\r^n \cup \r P^{n-1}$. $f$ then defines a 
projective (n-2)-plane $H_{infty} = \{ (x_1: \ldots :x_n) \in \r P^{n-1} | 
f(x_1, \ldots, x_n) = 0 \} \subset \r P^n$. This is in fact, the 
intersection of the 
projective hyperplanes $H_c = \overline{ \{ x \in \r^n | f(x) = c \} }$.

Now we are ready to state the genericity condition $f$ has to satisfy 
with respect to $M$. In next subsection we prove that

\begin{Th} \label{trans-condition}
If the linear functional $f$ is such that $H_{\infty}$ is transverse
to all the strata $X_1, \ldots, X_n$ at infinity, then the implications of 
the theorems A$'$ abd B$'$ hold for $f_{|M}$.
\end{Th}

\begin{Rem}
Given a submanifold $X \subset \r P^n$, almost any plane is 
transverse to $X$, same is true if we have a finite number of 
submanifolds. Hence for a generic functional $f$ the above 
transversality condition is satisfied.
\end{Rem}

\subsection{Thom's isotopy lemma and the proof of the main theorem}

There is a generalization of the theorem (A) of classical Morse theory to 
Whitney stratified space. It is the so-called \textit{Thom's first isotopy lemma}. 

\begin{Th}[Thom's First Isotopy Lemma]
Let $Y$ be a subset of a smooth
manifold $X$ with Withney stratification $\mathcal A$ and let $f: X \rightarrow P$ 
be a smooth map into another smooth manifold $P$, such that for each stratum 
$S \in \mathcal A, f_{|S}$ is a 
submersion and $f_{|\overline{S} \cap Y}$ is a proper map. Then $f: Y
\rightarrow P$ is locally trivial over $P$.
\end{Th}

The situation in the isotopy lemma is not exactly as in the theorem (A$'$). Because 
even though $\overline{M} \subset \r P^n$ is a compact Whitney stratified set,
$f : M \rightarrow \r$ can not be extended in a good way to the whole $\overline{M}$ to 
obtain a proper map on strata. Although with a simple trick we can convert the 
situation so that we can apply the isotopy lemma. 

Using the isotopy lemma, we prove the following theorem which immidiately 
implies theorem (A$'$).  

\begin{Th}[Main theorem]
Let $X$ be a compact smooth manifold and $D$ a domain in $\r ^k$ and 
$F:X \times D \rightarrow X$ 
a smooth map such that $F_z:X \rightarrow X$
is diffeomorphism for any $z \in \r ^k$.  
Suppose $M \subset X$ is a compact submanifold of $X$
of codimension $k$. $F_z(M)$ can be thought of as a 
family of submanifolds of $X$ parametrized by $\r ^k$.

Also suppose $Y \subset X$ is a Withney stratified subset of $X$
such that the following transversality 
condition holds:

$\forall z \in D$ we 
have $F_z(M)$ is transversal to all strata $S$ in $A$.

Then for all $z_1, z_2 \in D, F_{z_1}(M) \cap Y$ is homeomorphic to
 $F_{z_2}(M) \cap Y$ under an stratum 
preserving homeomorphism.
\end{Th}

\begin{Rem}
To obtain theorem (A$'$) from the above theorem, we take $X$ to be $\r P^n$, $k=1$, 
$M = \overline{ \{x \in \r^n | f(x) = 0 \} } \cong \r P^{n-1}$ and $D = [a, b]$. For 
$x = (x_0, \ldots, x_n)$ define 
$$ F(x, r) = (x_0, \ldots, x_{n-1}, x_n / r). $$

Then $F_z(M)$ is simply the closure of a level set of $f$, i.e. 
$$ F_z(M) = \overline{f^{-1}(z)}.$$

Hence the theorem implies that $\overline{f^{-1}(z)} \cap \overline{M}$ 
(for $a \leq z \leq b$) 
are homeomorphic under a stratum preserving map.
In particular, $f^{-1}(z) \cap M$ are homeomorphic.  
\end{Rem}

\begin{proof}
Consider $\hf: M\times D \rightarrow X\times \r ^k$ given by 
$\hat{F}(m,z) = (F(m,z), z)$. Obviously $\hat{F}$ is an embedding.
Now consider $Y \times D \subset X \times \r^k$ as a Whitney stratified 
space with product stratification $A \times \{\r\}$ (cf. ~{intersection-start}). 
Since by assumption
$F_z(M)$ is transversal to all strata in $Y$, so
$\hf$ is transversal to all strata in $Y \times D$. Thus inverse image of 
of starta in $Y \times D$ under $\hf$ gives a Whitney stratification 
for $\hf ^{-1}(Y \times D) \subset M\times D$ (cf. ~{intersection-of-strat}). 
Now consider projecton 
$t : M\times \r^k \rightarrow \r^k$ as a functon on $\hf^{-1}(Y \times
\r^k)$. Then since for any $z$, $F_z: M \rightarrow X$ is an embedding 
and $F_z(M)$ is transversal to all strata in $Y$ so $t$ is submersion 
restricted to any stratum of $\hf ^{-1}(Y \times D)$, 
and thus by Thom's First Isotopy Lemma, $t: \hf^{-1}(Y \times D)
\rightarrow D$ is a locally trivial fiberation. Notice 
$t ^{-1}(z) = \{ (m,z) | F_z(m) \in Y \} = F_z(M) \cap Y$. So for $z_1,
z_2 \in D, F_{z_1}(M) \cap Y$ is homeomorphic to $F_{z_2}(M) \cap Y
$ under a stratum preserving homeomorphism. 
\end{proof}

\section{Generalized Morse theory of Palais and Smale and an 
alternative proof of the theorems A$'$ and B$'$}

To be able to apply Morse theory to infinite dimensional spaces of
loops and geodesics, Palais and Smale developed a generalized Morse 
theory (cf. ~{Palais-Smale} and ~{Palais}). Since infinite 
dimensional manifolds are not compact, their theory can be also 
applied to the non-compact finite dimensional manifolds.

Following ~{Palais-Smale} (\textit{Generalized Morse theory}),
we give a brief account of their theory:

Let $M$ be a $C^2$-Riemannian manifold modeled on a separable Hilbert
space (hence it can be infinite dimensional). Let $f: M \rightarrow \r$ 
be a $C^2$ function. Assume that $f$ satisfies the following extra condition

\textit{
(C) If $\{ x_i \}$ is a sequence in $M$ on which $|f|$ is bounded and  
$\| \nabla f(x_i) \|$ converges to zero, then there is a critical 
point of $f$ in the closure of the set $\{ x_1, x_2, \ldots \}$. 
}

\begin{Th}
Let $M$ and $f$ satisfy the condition (C) above and assume that all the 
critical points of $f$ are non-degenerate. Then
\begin{enumerate}
\item For any real numbers $a < b$ there are only finitely many critical 
points of $f$ in $M_{a,b} = \{ x \in M | a < f(x) < b \}$, hence the 
critical values of $f$ are isolated.
\item Let $a$ and $b$ be regular values of $f$ and suppose that among the 
critical points of $f$ in $M_{a,b}$ there are $r$ having finite index.
Let the indices of these critical points be $\lambda_1, \ldots, \lambda_r$.
Then $M_{\leq b}$ has the homotopy type of $M_{\leq a}$ with $r$ cells
of dimensions $\lambda_1, \ldots, \lambda_r$ attached.
\end{enumerate}
\end{Th}

We now prove that in the situation of theorem ~{trans-condition}, 
the condition (C) is satisfied and hence Palais-Smale theory gives an
alternative proof of out theorems (A$'$) and (B$'$).

\begin{Prop}
Let $M$ be a (closed) algebraic submanifold of $\r^n$ of dimension $d$
and $\mathcal A = \{ X_0, X_1, \ldots X_r \} (X_0 = M)$ a finite 
Whitney stratification for $\overline{A} \subset \r P^n$.
Suppose $f$ is a linear functional such that $H_{infty} = 
\{ (x_1: \ldots :x_n) \in \r P^{n-1} | f(x_1, \ldots, x_n) = 0 \}$ 
is transverse to all the strata $X_1, \ldots, X_n$ at infinity.
Then the condition (C) holds for $M$ and $f$.
\end{Prop}
 
\begin{proof}
Let $\{ x_i \}$ be the sequence in the condition (C). If $\{ x_i \}$ is 
bounded there is nothing to prove. So assume it is not bounded. Without
loss of generality assume $x_i$ converges to $x \in \overline{M} \setminus M$
and $\lim_{i \rightarrow \infty} \| \nabla f_{|M} (x_i) \| = 0$. Now the
sequence of tangent planes $T_{y_i}M$ should have a convergent subsequence
in the Grassmanian $Gr(n, d)$ of $d$-planes in $\r^n$ as Grassmanian is compact.
Again without loss of generality assume that $T_{x_i}M$ converges
to $L \in Gr(n, d)$. Since $\| \nabla f_{|M}(x_i) \|$ goes to zero 
as $i$ goes to infinity, $Angle(\nabla f, T_{x_i}M)$, 
the angle between $\nabla f$ and $T_{x_i}M$, should also tend to
zero. But since $Angle(\nabla f, \cdot)$ is continuous we must have
$Angle(\nabla f, L) = 0$, that is $L \subset H_{\infty}$. By Whitney
condition (A) if $\lim_{i \rightarrow \infty} T_{x_i}M = L$ and 
$\lim_{i \rightarrow \infty} x_i = x$ then $T_x X \subset L$, where $X$ 
is the strata containing $x$. Hence $T_x X \subset H_{\infty}$. But this 
is a contradiction as $X$ and $H_{\infty}$ are transverse.
\end{proof}

\section{Morse theory for the orbits of a Lie group action}
In this section we reconsider the Morse theory of the section
\~ref{} for the case of orbits of Lie group actions.
We will give stronger versions (in the sense to be explained) 
of the isotopy lemma and our main theorem for this case as well as direct 
self-contained proofs.

Orbits of actions of Lie groups provide  interessting examples of Whitney
stratifications

\begin{Th}[cf. ~{}]
Let a Lie group $G$ acts smoothly on a manifold $X$ and let 
$Y \subset X$ be a closed invariant subset consisting of a 
finite number of orbits. Then the decomposition of $Y$ into 
orbits is a Whitney stratificatin of $Y$.
\end{Th}

Morse theory for stratified sets (Thom's first isotopy lemma) guaranttees
that the level sets of a Morse function on a stratified set are
homeomorphic under a stratum preserving homeomorphism, provided that 
certain transversality condition holds. Our generalization of 
isotopy lemma also asserts that intersections of a stratified set with
a moving family of submanifolds are homeomorphic as long as 
certain transversality condition holds.

But in general, we can not talk about this sections be diffeomorphic
under a  stratum preserving diffeomorphism as the Whitney 
stratification does not carry enough smooth straucture. The 
following famous example illustrates this

\begin{Ex}
Let $\gamma: \r \rightarrow (0, \infty)$ be a smooth 
non-constant function and let $Z \subset 
\r^3$ be defined by the equation
$$ xy(x+y)(x-\gamma(z)\.y)=0.$$
$Z$ is stratified by the $z$-axis which we call it $X$ and its 
complement (denoted by $Y$). This is in fact a Whitney stratification.
Adding $\r^4 \setminus A$ as a third stratum we obtain a  Whitney 
stratification of $\r^3$. Now consider the function $\pi_z$, projection 
on $z$ coordinate. The level sets of $\pi_z$ are planes parallel to 
the $xy$-plane. The intersection of a level set with the strata 
consists of four lines $x$-axis, $y$-axis, the line $x+y=0$ and 
the line $x-\gamma(z)\.y = 0$. 

It is obvious that the topological picture of the sections of 
the stratification remains the same as $z$ is changing, but now let us 
see if the picture is 
the same in the smooth category, i.e. if one can find a 
stratum preserving diffeomorphism 
between two different level sets ${\pi_z}^{-1}(c_1)$ and ${\pi_z}^{-1}(c_2)$.
This is in fact impossible 
because if there is such a diffeomorphism, its derivative at 
origin will be a linear transformation in $\r^2$ which leaves 
the three lines $x$-axis, $y$-axis and $x+y=0$ invariant and 
maps $x-\gamma(c_1)\.y = 0$ to $x-\gamma(c_2)\.y = 0$. But this is 
impossible since a linear transformation preserve the cross ratio of 
a collection of $4$ lines through origin in $\r^2$ while the 
cross ration of our four lines is cahnging according to $\gamma(z)$
(unless $\gamma$ is constant).  
\end{Ex}

In the proof of the first theorem of classical Morse theory 
(thoerem A) one uses the gradient vector field of the Morse function
$f$ and constructs a vector field such that $f$ increases with constant 
speed along trajectories of the vector field. The flow of this vector 
field then gives us the required diffeomorphism. To prove the isotopy lemma
, one basically follow the same idea but one needs the vector fields be
tangent to the strata as well so that the flow of the vector 
fields preserve the strata. But the strata of a Whitney stratification 
are not attached together in a smooth way hence we can only 
construct \textit{continuous} vector fields (rather than smooth) 
tangent to the stata and satisfying our required properties. Thus the 
flows of the vecor fields give us only homeomorphisms. 

But in the case the strata are the orbits of a Lie group action, there is 
a priori lots of smooth vector fields tangent to all the orbits, i.e. 
generating vector fields of the action. This enables us to prove stronger
version of our isotopy lemma (i.e. diffeomorphism insetad of
homeomorphism). 

The following can be considered as \textit{Thoms' first 
isotopy lemma for Lie group actions}

\begin{Prop}
Let a Lie group $G$ act smoothly on a manifold $X$ and let $Y 
\subset X$
be an invariant subset with compact closure. Also assume $f: 
X \rightarrow \r$ is a smooth function, suppose $[a, b] \subset \r$ such 
that $\forall c \in [a, b], f^{-1}(c)$ is transversal to all orbits in $Y$, 
then there exists vector field  $v$ with compact support on $X$ satisfying:
\begin{enumerate}
\item For any $y \in Y, v(y)$ is tangent to orbit of $y$.
\item Lie derivative of $f$ along $v$ is equal to $1$ at any point $y \in Y$.
\end{enumerate}
existence of $v$ then implies that for $c_1$ and $c_2 \in [a, b], 
f^{-1}(c_1) \cap Y$ 
is diffeomorphic to $f^{-1}(c_2) \cap Y$ under an orbit preserving 
diffeomorphism of whole space $X$.
\end{Prop}
\begin{proof}
Let $y \in \overline{Y}$ then,

Claim: There exist an open neighbourhood $U_y$ of $y$ and a 
smooth 
non sero vector field $\xi_y$ on $U$ such that $\xi_y$ at any $x$ in $U$ 
is tangent to orbits of $x$ and Lie derivative of $f$ along $\xi_y$ is 
non zero on $U$.

Proof of Claim: For any $\xi \in \mathfrak{g} = T_eG$, let 
$v_\xi$ be the 
left invariant vector field on $X$ generated by $\xi$ (i.e. $v_\xi(x)
= d \lambda (e) (\xi)$ where $\lambda : G \rightarrow X, \lambda (g) = g.x$). 
$v_\xi$ at $x$ is obviously tangent to orbit of $x$.
Since $\xi \mapsto v_\xi(y)$ is surjective , as a map from $\mathfrak{g}$ to 
$T_yO_y$, there exist $\xi \in \mathfrak{g}$ such that Lie derivative of
$f$ along $v_\xi$ at $y = L_{v_\xi}f(y) = df(d\lambda(e)(\xi))$ is 
non zero .
Let's fix $y$ and denote $v_\xi$ by $\xi_y$. By continuity we can find a 
nieghbourhood $U_y$ of $y$ such that $L_{\xi_y}f$ is non zero on $U_y$. 

Let $v_y = \xi_y / (L_{\xi_y}f)$, then obviously Lie derivative of 
$v_y$ on $U_y$ is 1. Now we use partition of unity to patch all $v_y$s 
together to get a compactly supported vector field $v$ on a 
neighbourhood of $Y$ such that Lie derivative of $f$ along $v$
is $1$, flow of $v$ then gives us the desired diffeomorphism 
between $f^{-1}(c_1) \cap Y$ and $f^{-1}(c_2) \cap Y$ (for $c_1$ and 
$c_2$ in $[a, b]$.
 
$\overline{Y} \subset X$ is compact . Take neighbourhood $W_y$ of $y$ 
such that $\overline{W_y} \subset U_y$, and take $W_{y_1}, \ldots, W_{y_n}$
among $W_y$s such that they cover $\overline{Y}$. Let $\eta_1, \ldots \eta_n
$ be partition of unity corresponding to open sets $\{ U_{y_1}, \ldots U_{y_n}
 \}$ such that $supp(\eta_i) \subset W_{y_i}$. each $\eta_i$ can be smoohtly 
extended to $X$ by letting $\eta_{y_i}$ to be zero outside $U_{y_i}$.
Let $v = \sum \eta_{y_i}v_{y_i}$ then $v$ is defined every where and 
$supp(v) \subset \cup \overline{W_{y_i}}$ which is compact. so $v$ is 
compactly supported.

\begin{eqnarray*}
L_{v}f(x) &=& df(x)(\sum_{i} \eta_{y_i}(x).v_{y_i}(x)) \cr 
&=& df(x)(\sum_{\{ i | x \in W_{y_i} \}} \eta_{y_i}(x).v_{y_i}(x) \cr
&=& \sum_{\{ i | x \in W_{y_i} \}} \eta_{y_i}(x).L_{v_{y_i}}(x) \cr
&=& \sum_{\{ i | x \in W_{y_i} \}} \eta_{y_i}(x) \cr
&=& 1 
\end{eqnarray*}
\end{proof}

Using the same idea that used was above before we can prove a stronger 
version of our main theorem in subsection ~{}

The following theorem is the main theorem we need and is basically a 
generalized version of the above propostion where instead of level set 
of function $f$ i.e. $f^{-1}(z)$ we consider a smooth family of
submanifolds $M_{z}$ parametrized by $z \in \r$. It would be useful to
generalize the situation and consider the family to be parametrized by
$\r^k$, this is given by a smooth map $F: X \times \r^k \rightarrow X$
where $F_{z}: M \rightarrow X$ 
is a diffeomorphism for each $z \in \r^k$ and we consider the smooth 
family $M_{z}$ to be $F_{z}(M)$. 

More precisely:

\begin{Th}[Main Theorem]
Let $X$ be a compact real manifold and 
$F:X \times\mathbb{R} ^k \rightarrow X$ 
smooth such that $F_z:X \rightarrow X$
is diffeomorphism for any $z \in \mathbb{R} ^k$.
Suppose $M \subset X$ is a compact submanifold of $X$
of codimension $k$. 
Also assume that the Lie group $G$ acts on $X$ smoothly and 
let $Y \subset X$ be a closed invariant subset, and 
finally assume that the following transversality
condition holds:

$\forall z \in D = I_1 \times I_2 \times \cdots \times
I_k \subset \mathbb R^{k}$ and $\forall y \in Y$, we
have $F_z(M)$ is transversal to $O_y$, (where $O_y =$ orbit of
$y$ and $I_i$s are closed intervals). 

then there exist vector fields $v_1, \ldots ,v_k$ on 
$X \times \mathbb R^{k}$ satisfying:
\begin{enumerate}
\item $v_i$ has compact support $(i = 1, \ldots ,k)$.
\item If $\widehat{F} (x, z) = (F(x, z), z)$ and
$\widehat{M} = \widehat{F}(M \times \mathbb R^{k})$
then $\forall (x,z) \in \widehat{M}$, $v_i(x, z) 
\in T_{(x,z)} \widehat{M}$, i.e. $v_i$ s are tangent 
to $\widehat{M}$, ($i = 1, \ldots ,k$).

\item For any $(x,z) \in X \times D, v_i(x,z) \in T_xO_x \times
\mathbb{R} ^k = T_{(x,z)}(O_x \times \mathbb{R} ^k)$, i.e.
$v_i$s are tangent to orbits of $G$ acting on $X \times D$
(where $G$ acts trivially on $\mathbb{R} ^k$).
\item (Main Property) $L_{v_i}t_j(x,z) = \delta_{ij}$, $\forall
(x,z) \in \widehat{M} \cap (Y \times D)$, where $t_j: X \times 
\mathbb{R} ^k \rightarrow \mathbb{R}$, $t_j(x, z) = z_j = j$-th 
component of $z$, and $L$ is for Lie derivative.
\end{enumerate}
In addition the existence of $v_i$s implies that all $F_z(M) \cap
Y$ 
are diffeomorphic for $z \in D$.
\end{Th}

\begin{proof}
Suppose $v_i$s are constructed satisfying $(1)$ to $(4)$. Let $\phi_i$ be
flow
of $v_i$, since $v_i$ are compactly supported $\phi_i ^t$ is
defined for all $t$. Note that $v_i$ is tangent to $\widehat{M}$,
so defines a vector field on $\widehat{M}$ and so $\phi_i ^t$ maps
$\widehat{M}$ to $\widehat{M}$.
Now consider $\psi :\widehat{M} \times \mathbb{R} ^k \rightarrow \widehat{M}$,
given by $\psi (\hat{x}, t_1, \ldots, t_k) = \phi_1 ^{t_1} (\phi_2 ^{t_2}(
\ldots \phi_k^{t_k}(\hat{x}) \ldots )$, for $\hat{x} \in \widehat{M}$ and
$(t_1, \ldots ,t_k) \in D$.

We show that $\psi_z: \widehat{M} \rightarrow \widehat{M}$ is a
diffeomorphism, where $z = (t_1, \ldots, t_k)$ which maps $M_{z_0} \cap
(Y \times D) $ onto $M_{z_0+z} \cap (Y \times D)$ for $z_0$ and $z_0 + z
\in D$. This is because of $(1)$ to $(4)$. Since $v_i$s are tangent to
orbits of
$G$ acting on $X \times \mathbb{R} ^k$ and $Y \times \mathbb{R} ^k$ is a union of
orbits so
$\phi_1^t, \ldots ,\phi_k^t$ map $Y \times \mathbb{R} ^k$ onto itself and so
does
$\psi_z$. Same reasoning shows that $\psi_z$ maps $\widehat{M}$ onto
itself.
Now Lie derivative of $t_i$ along $v_j$ is equal to $\delta_{ij}$,
this means that $t_i( \psi_z(\hat{x})) = t_i( \hat{x}) + z_j$,
where $z = (z_1, \ldots, z_n)$. So $\psi_z$ maps $M_{z_0} \cap
\widehat{M} = \{ \hat{x} \in \widehat{M} : t_i(\hat{x}) = z_i \}$ onto 
$M_{z_0+z} \cap \widehat{M} =
\{ \hat{x} \in \widehat{M} : t_i(\hat{x}) = z_{0,i} + z_i \}$

Now we show how to construct $v_i$s satisfying $(1)$ to $(4)$.
Recall $\widehat{F}(x, z) = (F(x, z), z)$. Let $ \chi_i(x, z) = 
d\widehat{F}_{(x,z)}(0, e_i)$ be the velocity vector fields, where
$(0, e_i) \in T_xX \times \mathbb{R}
^k$ and $e_i$ is the $i$-th standard basis element in $\mathbb{R} ^k$.
Then since $\widehat{F}(M \times \mathbb{R} ^k) = \widehat{M}, \chi_i$
are tangent to $\widehat{M}$ and so $\chi_i$ restricted to $\widehat{M}$
defines a vector field on it. 

Now take $(x, z) \in \widehat{M} \cap (Y \times D)$, let $(x, z) = 
\widehat{F} (p, z)$ where $p \in M$ and let $\gamma_1, \ldots,
\gamma_m$ be vector fields on a neighbourhood $U$ of $p$ in $M$ such
that $\gamma_1, \ldots, \gamma_m$ generate $T_yM$, $ \forall y \in U,
(m =dim(M))$,  and let $\nu_i (F(y, z), z) = d\widehat{F}_{(y,z)} 
(\gamma_i(y), 0)$ be the image of $\gamma_i$ under $\widehat{F}$.
Since $\widehat{F}_z$ is embedding $\nu_1 , \ldots, \nu_m$ give a 
basis for $T_{\widehat{F}(y,z)} (\widehat{F}_z(M)), \forall y \in U$.

Now since $O_x$ is transversal to $F_z(M), \forall x \in Y, \forall
z \in D$, one can find $\xi_1, \ldots, \xi_k \in \mathfrak{g} =$ Lie
algebra of $G$ such that the invariant vector fields on $X \times 
\mathbb{R}
^k$ generated by $\xi_1, \ldots, \xi_k$ (which we again denote them by
$\xi_1, \ldots, \xi_k$) together with $\nu_1, \ldots \nu_m$ and $
\partial / \partial{t_1}, \ldots, \partial / \partial{t_k}$ form a
basis for $T_{(x,z)} (X \times \mathbb{R} ^k)$, where $(\partial
/ \partial{t_1}, \ldots, \partial / \partial{t_k})$ is the standard
basis for $\mathbb{R} ^k$.
Thus there exists a neighbourhood $U$ of $(x, z)$ in $X \times \mathbb{R}
^k$ such that $\xi_1(y, \zeta), \ldots, \xi_k(y, \zeta)$ together with
$\nu_1(y, \zeta), \ldots \nu_m(y, \zeta)$ and $
\partial / \partial{t_1}, \ldots, \partial / \partial{t_k}$ form a
basis for $T_{(y, \zeta)} (X \times \mathbb{R} ^k)$, for all $(y, \zeta)
\in U$. Now let vector fields $V_{(x, z)} ^i$ be the projection of
$\chi_i(y,
z)$ on the space generated by $\xi_1, \ldots , \xi_k$ and $\partial
/ \partial{t_1}, \ldots, \partial / \partial{t_k}$, having defenition
of $\chi_i$ in mind, one can write then $V_{(x, z)}^i(y, \zeta) =
\alpha^1 (y,
\zeta) \xi^1(y) + \cdots + \alpha^k (y, \zeta) \xi^k(y) + \partial /
\partial{t_i}$.

$V_{(x, z)}^i$ are vector fields on $U$ such that:
\begin{enumerate}
\item $V_{(x, z)}^i (y, \zeta) \in T_yO_y \times \mathbb{R} ^k$
\item $L_{t_j}V_{(x, z)}^i (y, \zeta) = \delta_{ij}$ because 
$L_{t_j}\chi_i = \delta_{ij}$ and $L_{t_j}\nu_i = \delta_{ij}$
\item $V_{(x, z)}^i (y, \zeta) \in T_{(y, \zeta)} \widehat{M}$ for 
$(y, \zeta) \in \widehat{M}$ because $\chi_i$ and $\nu_i$ are tangent 
to $\widehat{M}$.
\end{enumerate}

So for any $(x,z) \in \widehat{M} \cap (Y \times D)$ we get a 
neighborhood $U_(x,z) \in X \times \mathbb{R} ^k$ and vector fields
$V_{(x,z)}^i (i = 1, \ldots, k)$ defined on $U_{(x,z)}$. Notice 
that $\widehat{F}(M \times D)$ is compact, so $\widehat{M} \cap (Y \times
D ( \subset \widehat{F}(M \times D))$ is compact, Y being closed.

Take $W_{(x,z)} \subset W'_{(x,z)} \subset U_{(x,z)}$ such that 
$\overline{W} \subset W'$ and $\overline{W'} \subset U$, with $\overline{W}$ and
$\overline{W'}$ compact. Since $\widehat{M} \cap (Y \times D)$ is compact 
finitely many of $W$s cover $\widehat{M} \cap (Y \times D)$, denote them 
by $W_i, (i = 1, \ldots, r)$. Let $\eta_1, \ldots, \eta_r$ and $\eta$ be
partition of unity for $X \times \mathbb{R} ^k$ corresponding to 
open cover $W'_1, \ldots, W'_r$ and $(\bigcup_{i=1}^{r} \overline{W_i})^c$,
(since $\bigcup \overline{W_i} \subset \bigcup W'_i, W'_1, \ldots, W'_r$ and $
(\bigcup \overline{W_i})^c$ cover the whole space $X \times \mathbb{R} ^k$).
Denote vector fields $V_{(x,z)}^j$ corresponding to open set $U_i$ by 
$V_i^j, (j = 1, \ldots, k)$. Since $ supp(V_i^j) \subset U_i$ and $supp
(\eta_i) \subset W'_i$ and $\overline{W'_i} \subset U_i$ then
$\eta_i(x,z)V_i^j(x,z)$ are globally defined and smooth everywhere on 
$X \times \mathbb{R} ^k$ with $supp(\eta_i V_i^j) \subset \bigcup
W'_i$. Define $v_j(y, z) = \sum_{i=1}^{r} \eta_i(y,z) V_i^j(y,z)$,
then $v_j$ are smooth everywhere and $supp(v_j) \subset \bigcup W'_i$.
For $(y,z) \in \bigcup W_i, L_{t_k}v_j(y,z) = \sum_{i: (y,z) \in W'_i}
\eta_i (y,z) L_{t_k}V_i^j(y,z) = \sum_{i: (y,z) \in W'_i} \eta_i(y,z)
 \delta_{kj} = \delta_{kj}$, because $\eta_i$s are partition of unity 
with respect to $W'_i$ and $(y,z) \notin (\bigcup \overline{W_i})^c$. 

So for $(y,z) \in \bigcup W_i = $ an open set in $\widehat{M} \cap (Y
\times D)$ we have $L_{t_i}v_j(y,z) = \delta_{ij}$ and since $supp(v_j)
\subset \bigcup W'_i \subset \overline{W_i}$, which is compact, $v_j$s
are compactly supported. finally because of $(1)$ to $(3)$ for $V_i^j$s
we have:
\begin{enumerate}
\item $v_j(y,z) \in T_yO_y \times \mathbb{R} ^k$.
\item $v_j(y,z) \in T_{(y,z)}\widehat{M}$.
\\
and we just proved: 
\item $L_{t_i}v_j(y,z) = \delta_{ij}$ for $(y,z) \in$ a neighbourhood
of
$\widehat{M} \cap (Y \times D)$.
\end{enumerate}

\end{proof}

\chapter{Euler Characteristic of Sections of Orbits of Algebraic
Groups}

\section{Algebraic torus actions}

Notation: We denote the direct product of $n$ copies of
multiplicative group of complex numbers by $\tor$. Let $ x = (x_1,
\ldots, x_n)$ be an element of this group and let $k = (k_1,
\ldots, k_n)$ be an $n$-tuple of integers. We denote the monomial
$x_1^{k_1}\ldots x_n^{k_n}$ simply by $x^k$.

As we mentioned in the beginning, there is a beautiful theorem due
to D.N. Bernstien (cf. {Bernstein} and {Askold})
regarding the topology of hypersurfaces defined  in $\tor$. Let us
recall the theorem.

\begin{Def}
Let $f(x_1, \ldots, x_n)$ be a Laurant polynomial (i.e. negative
powers of $x_i$'s are allowed) with complex coefficients. To each
monomial $c_k x^k$ we can assign a point $k = (k_1, \ldots, k_n)$
in $\z ^n$. The \textit{Newton polyhedron} $\Delta$ of $f$ is
defined to be the convex hull of the points $k$ in $\z ^n$
corresponding to monomials of $f$.
\end{Def}

Newton polyhedron is a very nice combinatorial invariant of  $f$
which contains a lot of information about geometry and topology of
geometric objects defined by polynomials. It can be thought as the
generalization of degree of a one variable polynomial. The philosophy is 
that in the generic cases, one can describe all the discrete topological and 
geometric invariants of $f$ in terms of its Newton polyhedra.

For each polyhedron $\Delta \subset \r ^n$ and a covevtor $\xi$ we
define the polyhedron $\Delta ^\xi$ to be the face of the
polyhedron $\Delta$ on which $\xi$ attains a minimum (in
particular, $\Delta ^0$ is $\Delta$). If $\Delta$ is Newton
polyhedron of a polynomial $f$, then the sum of monomials
corresponding to points in $\Delta ^\xi$ is denoted by $f^\xi$.

\begin{Def}
We say that a system of Laurant polynomials $f_1, \ldots, f_k$ is
\textit{non-singular} for their Newton polyhedra, if for any
covector $\xi \in \r^n$ the following condition holds: for any
solution $z$ of the system $f_1^\xi (x) = \ldots = f_k^\xi (x) =
0$ in $\tor$, differentials $df_i ^\xi$ ($i=1, \ldots ,k$) are
linearly independent.
\end{Def}

One can show that non-singularity condition is a generic
condition, i.e. in the space of all systems of functions $f_1,
\ldots, f_k$ with fixed Newton polyhedra $\Delta_1, \ldots,
\Delta_k$, non-singularity condition holds for every system except
for a set of measure zero.

\begin{Th}[Bernstein ~{Bernstein}]
Let $f(x_1, \ldots , x_n)$ be a Laurant polynomial which is
non-singular for its Newton polyhedron $\Delta$. Let $X = \{
f(x_1, \ldots ,x_n) = 0 \}$. Then $\chi(X) = n! \cdot
Vol(\Delta)$.
\end{Th}

One can generalize the above theorem to complete intersections of 
hypersurfaces. We need a bit of notation. Let $\Delta_1, 
\ldots, \Delta_n$ be n-polyhedra in $\r^n$. Let $V(\Delta_1, \ldots, \Delta_n)$ 
denotes the mixed volume 
\footnote{For definition of mixed volume look at for example ~{}}
of these polyhedra. Now let $F(x_1, \ldots, x_k)$ be the Taylor series of 
an analytic function in $k$ variables at the point $0$. We wish to define 
the number $F(\Delta_1, \ldots, \Delta_k)$. If $F$ is a monomial of degree
$n, F(x_1, \ldots, x_n) = {x_1}^{n_1} \cdots {x_k}^{n_k}$, we put
$$F(\Delta_1, \ldots, \Delta_k) = n! V(\Delta_1, \ldots, \Delta_1, 
\ldots, \Delta_k, \ldots, \Delta_k).$$
where each $\Delta_i$ is repeated $n_i$ times in the mixed volume.
One extends the definition to a homogeneous polynomial $F$ of degree $n$
by linearity and for a power series $F$ one defines
$$F(\Delta_1, \ldots, \Delta_k) = n! F_n(\Delta_1, \ldots, \Delta_k).$$ 
where $F_n$ is the homogeneous part of $F$ of degree $n$.
We then have
\begin{Th}[Bernstein ~{Bernstein}]
Let $X$ be the variety defined oin $\tor$ by a non-degenerate system of 
equations $f_1 = cdots = f_k = 0$ with Newton polyhedra 
$\Delta_1, \ldots, \Delta_k$. Then 
$$\chi(X) = \Pi_{i=1}^{k}\Delta_i (1+\Delta_i)^{-1}.$$
\end{Th}
 
Bernstein's theorem can be put in a more fancy way as follow: let $p(x) =
\sum_{\alpha \in S}c_\alpha x^\alpha$ with $S = \{ \alpha_1,
\ldots, \alpha_N \} \subset \mathbb{Z} ^N$ be a Laurant
polynomial. As before we denote its Newton polyhedron by $\Delta$.
Let $\pi : \tor \rightarrow \Tor \subset GL(n, \c)$ be a finite
dimensional representation of algebraic torus $\tor$ given by
$\pi(x_1,\ldots, x_n) = diag( x_{1}^{\alpha_1}, \ldots,
x_{n}^{\alpha_n})$ and consider the hyperplane  $L = \{ x \in \c^N
| f(x) = 0 \}$ where $f(x_1, \ldots, x_n) = \sum_{i=1}^{n}
\alpha_i x_i$, then if $\alpha_i$ are generic enough, the Euler
characteristic of hyperplane section $L \cap
\pi(\mathbb{C}^{*^n})$
 is equal
to $n! Vol(\Delta)$.

On the other hand there is a remarkable formula due to
Kushnierenko ({Kushnierenko}) giving the number of solutions
of a system of Laurant polynomials.

\begin{Th}[Kouchnierenko]
Let $f_1, \ldots, f_n$ be $n$ Laurant polynomials in $n$ variables
with the same Newton polyhedron $\Delta$. If the coefficients of
the polynomials are generic enough, then the number of solutions of
the system $f_1(x)= \ldots =f_n(x)=0$ is equal to $n!
Vol(\Delta)$.
\end{Th}

There is a generalization of Kouchnierenko's theorem to any
system of generic polynomials (not necessarily with the same 
Newton polyhedra) due to D. Bernstein

\begin{Th}[Bernstein]
Let $f_1, \ldots, f_n$ be $n$ Laurant polynomials in $n$ variables
with Newton polyhedra $\Delta_1, \ldots, \Delta_n$. If the coefficients of
the polynomials are generic enough, then the number of solutions of
the system $f_1(x)= \ldots =f_n(x)=0$ is equal to $n!
V(\Delta_1, \ldots, \Delta_n)$.
\end{Th}

Kouchnierenko's theorem can also be formulated using 
representations of $\tor$. Let
us recall notion of degree of a subvariety $V$ of $\c ^n$ (or $\c
P^n$): if $dim(V)$ is $d$ then the number of intersections of $V$
and a generic plane of codimension $d$ is constant called
\textit{degree} of the subvariety $V$. Now let $\pi : \tor
\rightarrow \Tor \subset GL(n, \c)$ be the representation of
$\tor$ defined by monomials corresponding to the the points in
Newton polyhedron of $f_1, \ldots, f_n$. $\tor$ gets embedded in
$\c^N$ via $\pi$. Hence

\begin{Th}
Let us $\tor$ acts on $\c^N$ via a linear representation $\pi$.
Let $\Delta$ denotes the convex hull of weights of this
representation. Then degree of a generic orbit of this action is
equal to $n! Vol(\Delta)$.
\end{Th}

\section{Generalization of the Kouchnierenko's theorem
and the spherical varieties}.

\subsection{Generalization of Kouchnierenko's theorem to the representations 
of reductive groups}
Let $G$ be a connected $n$-dimensional complex reductive group and
let $\pi$ be an $N$-dimensional holomorphic representation of this
group. Consider the systems $f_1(x)  = f_2(x) = \cdots = f_n(x) =
0$, $x \in G$, where each $f_i$ is a linear combination of matrix
entries of the representation $\pi$. All the systems that lie
outside of a certain algebraic hypersurface in the space of such
systems, have the same number of roots which we denote by
$N(\pi)$.

Let $\t$ and $\t ^*$ be the Lie algebra and dual Lie algebra of a
maximal torus $T$ in $G$.

\begin{Def}
The Newton polyhedron $\Delta$ of the representation $\pi$ is
defined to be the convex hull of its weights.
\end{Def}

One can construct a homogeneous polynomial $V$ of degree $n$ on
the cone of all convex bodies in $\t ^*$, i.e. $V$ is the
restriction to the diagonal of an $n$-linear function on the cone
of convex bodies, such that
\begin{Th}[Kazarnovski(1986)]
$N(\pi)$ is equal to $n! \cdot V(\Delta)$. Or in other words, the
degree of the subvariety $\pi(G)$ of $M(N,\c)$, vector space of
all $N \times N$ matrices, is equal to $n! \cdot V(\Delta)$.
\end{Th}

For the construction of $V$ and the proof of this theorem look at
~{Kazarnovskii}.

\subsection{Further generalization of the Kouchnierenko's theorem:
spherical varieties}.

When we have a $N$-dimensional representation $\pi$ of $G$, we can
think of $G$ acting on the space of matrices $M(N, \c)$ via left
multiplication, or we can let $G \times G$ acts on $M(N, \c)$ by
left-right multiplication ($(g,h) \cdot m = g \cdot m \cdot
h^{-1}$). In both cases $\pi(G)$ would be the orbit of identity.
Hence Kazarnovskii's theorem gives a combinatorial formula for the
degree of this orbit.

It is interesting to see if there is similar combinatorial formula
for a bigger class of orbits of reductive groups. Brion (cf.
~{Brion}) has discovered such a formula for the so-called
\textit{Spherical varieties}. To give statement of his theorem we
need to introduce some notations and definitions:

\begin{Def}
Let $G$ be connected reductive group. A homogeneous space $G/H$ is
called spherical if a Borel subgroup of $G$ has a dense orbit in
$G/H$. Similarly a $G$-variety $X$ (i.e. a variety together with
an algebraic action of $G$) is called spherical if a Borel
subgroup of $G$ has a dense orbit.
\end{Def}

Let us see how Spherical varieties generalize some of the 
previously considered, interseting examples of group actions.
In all the examples $G$ is a complex reductive group.

\begin{Ex}
Let $G = \tor$. Then the Borel subgroup of $G$ is $G$ itself. Hence
$G$ is a spherical homogeneous space with natural left action of $G$
on itself, and spherical $G$-variety are just \textit{toric varieties}
by definition. 
\end{Ex}
 
\begin{Ex}
Supose $\sigma$ is an involution of $G$ (that is an automorphism of
order two of $G$). Let $H$ be the subgroup of all elements of $G$ fixed
by $\sigma$. The homogeneous space $G/H$ is called a \textit{symmetric
variety}. Origin of this notion, perhaps, is from differential geometry.
It is a manifold with high number of symmetries and hence interesting 
object to be studied. 
It is well-known that every symmetric variety is a spherical homogeneous 
space.
\end{Ex}

\begin{Ex}
Let $G \times G$ act on $G$ by left-right multiplication. From the 
so-called \textit{Bruhat decomposition}, one knows that $B \times B$
(which is a Borel subgroup of $G \times G$) has an open dense orbit in $G$. Hence
$G$ is a $G \times G$-spherical homogenuous space. This example is, in fact, what 
Kazarnovskii considered.  
\end{Ex}


Let$V$ be the representation space of a connected reductive
complex group $G$. Let $K$ be a maximal compact subgroup of $G$.
Action of $G$ on $V$ naturally induces an action on the projective
space $\p(V)$. Choose $K$-invariant K\"ahler structure on the
$\p(V)$, which is in particular, a symplectic manifold (the
symplectic form $\omega$ being the being the imaginary part of the
K\"ahler form). Define a map $\mu$ from $X$ to $\k ^{*}$ (the dual
of the Lie algebra $\k$ of $K$) by $\mu(x)(A) = (\tilde{x} \cdot
A\tilde{x})(\tilde{x} \cdot \tilde{x}) ^{-1}$ where $x \in \p(V)$;
$\tilde{x}$ is a representation of $x$ in $V$; $A \in \k$. It is
easy to see that $\mu$ is $K$-invariant and its differential
$d\mu$ satisfies $d\mu_x(\xi)(A) = \omega_x(\xi, A_x)$ for every
$x \in X, \xi \in T_x\p(V), A \in \k$. This means that $\mu$ is a
\textit{moment map} for the symplectic action of $K$ on $\p(V)$
(cf. ~{GS}).

If $X$ is any closed smooth algebraic subvariety of $\p(V)$ which
is $G$-stable, then $X$ inherits a $K$-invariant K\"ahler
structure, and the restriction of $\mu$ to $X$ is still a moment
map for the induced symplectic structure on $X$.

The image $\mu(X)$ has a nice convexity property. As $\mu(X)$ is a
$K$-stable subset of $\k^{*}$ ($K$ acting on $\k^{*}$ via
coadjoint representation), it is described by its intersection
with a fundamental domain $C$ of $K$ acting on $\k ^{*}$. We can
choose $C$ to be a Weyl chamber in the dual $\t ^{*}$ of Lie
algebra of a maximal torus of $K$.

\begin{Th}[~{Kirwan}, ~{GS}]
The intersection $\Delta = \mu(X) \cap C$ is a convex polyhedron
with rational vertices with respect to the weight lattice.
\end{Th}

$\Delta$ is usually called the \textit{moment polyhedron}. For the
case of action of algebraic torus $\tor$, the moment polyhedron is
same as the Newton polyhedron. One can think of moment polyhedron
as Newton polyhedron for spherical varieties.

 Let $R$ be a root system of $G$ (with respect
to the maximal torus associated to $\t$), $R^+$ the set of
positive roots defined by the choice of $C$, and $E$ the set of
all positive roots which are orthogonal to $\Delta$. We denote by
$\rho$ half the sum of all the positive roots.

The following theorem due to Brion (cf. ~{Brion}) give the
degree of $X$ as the integral of a certain function on the moment
polyhedron.

\begin{Th}[Brion (19??)]
The degree of an $n$-dimensional spherical subvariety $X$ of
$\p(V)$ is equal to
$$ n! \int_{\Delta} \prod_{\alpha \in R^{+} \setminus E} (\gamma ,
\alpha) / (\rho, \alpha) d\gamma$$.
\end{Th}

Following a suggestion by A.G. Khovanskii, A. Okounkov has define
a bigger polyhedron $\tilde{\Delta}$ over the moment polyhedron
$\Delta$ (i.e. it projects onto $\Delta$) such that the above
formula becomes $deg(X) = n! Vol(\tilde{\Delta})$ (cf.
~{Okounkov}).

Let $G$ be the torus $\tor$. Since $G$ is abelian it is equal to
its Borel subgroup and hence all homogeneous spaces of $G$ are
spherical.  A spherical $G$-variety then is a variety which has a
dense $G$ orbit. So If $G$ is the torus, spherical $G$-varieties
are exactly \textit{toric varities}. Consider a diagonal
representations of torus $G$ on $V = \c^N$ with weights $\omega_1,
\cdots \omega_N$. As before $G$ acts on $\p(V)$, let $X$ be the
closure of the orbit of $(1: \cdots :1) \in \p(V)$. One can easily
see that the moment polyhedron $\Delta$ in this case is exactly
the Nweton polyhedron, i.e. the convex hull of $\omega_1, \cdots
omega_N$ and also the above formula for degree of $X$ reduces to
$n! Vol(\Delta)$ yielding the Kouchnierenko's theorem.

One can also obtain Kazarnovskii's result as a special case of the
formula for degree of spherical varieties. For this, one considers
action of $G \times G$ on $G$ via multiplication from left and
right. Let $B$ denote a Borel subgroup of $G$. From a well-known
result in algebraic groups (i.e. Bruhat decomposition), $B \times
B$ which is a Borel subgroup of $G \times G$ has a dense orbit in
$G$ (the big Bruhat cell) and hence $G$ is a $G \times
G$-spherical orbit. Let $\pi$ be an $N$-dimensional faithful
representation of $G$. Then $G \times G$ can act on $M(N, \c)$ by
multiplication from left and right and $\pi(G)$ is the orbit of
identity. This action obviously induces an action of $G \times G$
on the projective space $\p(M(N, \c))$. Let $X$ be the closure in
$\p(M(N, \c))$ of the orbit of identity (that is image of $\pi(G)$
in the projective space). Then $X$ is a spherical variety since
$\pi(G)$ is a spherical $G \times G$ homogeneous space. If
$\omega_1, \cdots, \omega_t$ are weights of the representation
$\pi$, one can show that the moment polyhedron $\Delta$ is the
convex hull of the points $(\omega_i, -\omega_i), i=1, \cdots, t$
living in $\t^* \oplus \t^*$, that is the dual of the Lie algebra
of a maximal torus of $G \times G$ (here $\t$ is the Lie algebra
of a maximal torus of $G$). If we project $\Delta$ on $\t^*$ we
obtain $\Delta'$ the convex hull of weights $\omega_1, \cdots,
\omega_t$. This $\Delta'$ is the polyhedron appeared in the
Kazarnovskii's formula. In fact, if one rewrites the Brion's
formula as an integral over $\Delta'$ rather than $\Delta$, one
will recover Kazarnovskii's formula.

\section{Euler characteristic of hyperplane sections,the number of 
critical points and degree of a subvariety.}

\subsection{Euler characteristic and the number of critical points}
Let us recall the theorem
{Eulerchar-number-of-critical-points}, relating Euler
characteristic of a generic hyperplane section and the number of
critical points of a linear functional. Let $f$ be a linear
functional on $\c ^N$ and $Y \subset \c^N$ some submanifold, We
denote by $\mu(f, Y)$, the number of critical points of $f_{|_Y}$.

\begin{Th}
Let $X$ be a smooth closed algebraic subset of $\c^N$ of dimension
$d$ such that $\bar{X} \subset \c P^N$ has a Whitney
stratification with finite number of algebraic strata. Let $f$ be
a generic linear functional $f$ on $\c^N$ and $c$ a generic
complex number. We then have
$$ \chi(\{x \in X| f(x) = c \}) = \chi(X) + (-1)^{(d+1)} \cdot
\mu(f, X). $$
\end{Th}

Now suppose a complex Lie group $G$ acts linearly on $\c^N$, the
action obviously extends to $\c P^N$ by letting $G$ act trivially
on the last homogeneous coordinate. Let $X$ be a smooth closed
orbit of $G$ in $\c^N$ such that $\bar{X} \subset \c P^N$ consists
of a finite number of $G$-orbits. Then it is well-known that
$\bar{X}$ admits a Whitney stratification with $G$-invariant
strata (i.e. each stratum is union of orbits). So we can apply the
above theorem to get

\begin{Th}
Let a complex Lie group $G$ act linearly on $\c ^N$ (and hence on
$\c P^N$). Suppose $X \subset \c^N$ is a closed $G$-orbit of
dimension $d$ whose closure in $\c P^N$ consists of a finite
number of orbits. Then for a generic linear functional $f$ on
$\c^N$ and a generic complex number $c$ we have
$$ \chi(\{x \in X| f(x) = c \}) = \chi(X) + (-1)^{(d+1)} \cdot
\mu(f, X). $$
\end{Th}

Let us consider an interesting special case of this theorem. It
turns out that in this case the conditions of the above theorem
are easy to verify. Let $\pi : G \rightarrow GL(N, \c)$ be a
faithful representation of a complex connected reductive group
$G$. One can define a linear action of $G \times G$ on $M(N, \c)$,
the vector space of matrices, by $(g, h)\cdot m = \pi(g) \cdot m
\cdot \pi(h)^{-1}$. Then $G\times G$-orbit of identity will be
simply $\pi(G)$.

The following theorem tells us what $chi$ of $\pi(G)$ is.

\begin{Prop}
If $G$ is a complex connected reductive group then Euler
characteristic of $G$ is equal to zero.
\end{Prop}

\begin{proof}
By Iwasawa decomposition $G$ has the same homotopy type as its
maximal compact subgroup $K$. In particular $\chi(G) = \chi(K)$.
Any non-zero left invariant vector field on the compact Lie group
$K$ gives a everywhere non-zero vector field on $K$ and hence by
Poincare-Hopf theorem $\chi(K) = 0$.
\end{proof}

It is an easy lemma to prove that if one has a homomorphism from
an algebraic groups $G_1$ to another algebraic group $G_2$, image
of $G_1$ in $G_2$ is always a closed subset of $G_2$, hence in our
case $\pi(G)$ is always closed subset of $GL(N, \c)$, but it may
not be a closed subset of $M(N, \c)$. Next propositions, provide a
nice criterion to check when $\pi(G)$ is a closed subset of $M(N,
\c)$, namely $\pi(G)$ is closed in $M(N, \c)$ iff origin belongs
to the interior of the convex hull of weights of the
representation.

To prove the general case, we need the special case of torus
\begin{Prop}
Let $\pi: \tor \rightarrow GL(N, \c)$ be a representation of an
algebraic torus with weights $\omega_1 ,\cdots ,\omega_N \in
\z^n$. Then $\pi(\tor)$ is a closed subset of $M(N, \c)$ iff
origin belongs to the interior of the convex hull of the weights
$\omega_1 ,\cdots , \omega_N$.
\end{Prop}
\begin{proof}
See {Vinberg-Popov}.
\end{proof}

\begin{Prop}
Let $\pi: G \rightarrow GL(N, \c)$ be a representation of a
reductive group $G$. Then $\pi(G)$ is a closed subset of $M(N,
\c)$ iff $\pi(T)$ is closed in $M(N, \c)$ that is iff origin
belongs to the interior of the convex hull of the weights of the
representation $\pi$.
\end{Prop}

\begin{proof}
Let $T$ be a maximal torus of $G$ and $B$ a Borel subgroup
containing $T$, let B=TU be a Levi decomposition and let K be a
maximal compact subgroup of $G$ such that $G=TUK$. We may assume
that $\pi(B)$ is contained in the usual upper triangular Borel
subgroup of $GL(N, \c)$ and $\pi(T)$ is contained in the diagonal
matrices, and $\pi(U)$ is contained in the upper triangular
matrices with $1$'s down the diagonal. Suppose
$g_n=\pi(t_n)\pi(u_n)\pi(k_n)$ is a sequence in $\pi(G)$ that
converges to $g$ in $M(N,\c)$. with $t_n \in T, u_n \in U$ and
$k_n \in K$. We want to show that $g$ is in $\pi(G)$. Since
$\pi(G)$ is closed in $GL(N, \c)$, as was explained before, it
suffices to show that $g$ is in $GL(N,\c)$. Suppose not. Since
$\pi(K)$ is compact, we may assume that $\pi(k_n)$ converges to an
element $k$ in $GL(N,\c)$. We may thus assume that $g_n = \pi(t_n)
\pi(u_n)$. This is a sequence in the standard Borel subgroup of
$GL(N, \c)$ converging to an element $g$ not in $GL(N,\c)$. This
implies that $\pi(t_n)$ converges to an element not in $\pi(T)$,
contradicting the assumption that $\pi(T)$ is closed in $M(N,\c)$.
\end{proof}

As was mentioned in previous section, $G$ under left-right action
of $G \times G$ is an example of so-called spherical homogeneous
space, that is $B \times B$, which is a Borel subgroup of $G
\times G$, has a dense orbit in $G$ (Bruhat decomposition). Of
course this equally applies to $\pi(G)$ i.e. $\pi(G)$ is a
spherical $G \times G$ homogeneous space. It is a well-known
result in the theory of spherical varieties that any
compactification of a spherical homogeneous space, consists of a
finite number of orbits. In our case it means that closure of
$\pi(G)$ in $\c P^N$ consists of a finite number of $G \times G$
orbits.

Now we restate the theorem about the Euler characteristic of
generic hyperplane sections in the special case of a reductive
group embedded in vector space of matrices.

\begin{Th}
Let $\pi: G \rightarrow GL(N, \c)$ be a faithful representation of
a $d$-dimensional complex connected reductive group $G$. Suppose
$\pi(G)$ is closed in $M(N, \c)$ that is origin belongs to the
convex hull of weights of $\pi$. Then for $f$ a generic linear
functional on $M(N, \c)$ and $c$ a generic complex number we have
$$ \chi(\{x \in \pi(G)| f(x) = c \}) = (-1)^{(d+1)} \cdot \mu(f,
\pi(G)).$$
\end{Th}

\subsection{Proof of the Bernstein's theorem}

Suppose the torus $\tor$ is acting on $\c^N$ via a diagonal
faithful representation $\pi$. We can view $\pi$ as a homomorphism
$\pi : \tor \rightarrow \Tor \subset GL(N, \c)$, as the image of
$\pi$ lies in the subspace of diagonal matrices. so we can think
of the torus $\tor$ as embedded in a bigger torus $\Tor$ via the
homomorphism $\pi$.

As was mentioned before in section{}, one can reformulate the
Bernstein's theorem regarding Euler characteristic of a generic
hypersurface in $\tor$ as Euler characteristic of a generic
hyperplane section of the torus $\tor$ embedded in $\c^N$ is equal
to $(-1)^{(n+1)}$ times its degree as a subvariety of $\c^N$.

First of all we claim that we can assume that the image of $\tor$
in $\c^N$ is closed, or equivalently the origin is in the interior
of the Newton polyhedron of $\pi$. If not, we shift the Newton
polyhedron so that origin lies inside the polyhedron. Shifting the
polyhedron corresponds to multiplying all the monomials by a fixed
monomial. Since we are considering a hypersurface in $\tor$,
multiplying everything by a fixed monomial does not change the
hypersurface in $\tor$, i.e. without loss of generality we can
assume origin is inside the Newton polyhedron. As we mentioned in
the previous section this means that the image of $\tor$ in $\c^N$
is closed. We are now in the position to use theorem {}, to
get that Euler characteristic of a hyperplane section is equal to
$-1^(n+1) \cdot \mu(f, \pi(\tor))$, where $f$ is a generic linear
function defining the hyperplane in $\c^N$. So the only thing
remains to complete the proof of the Bernstein's theorem is to
show that number of critical points of a generic linear functional
restricted to the torus is equal to its degree

\begin{Prop}
Let $\pi: \tor \rightarrow \Tor \subset GL(N, \c)$ be a diagonal
faithful representation of $\tor$. Let $f$ be a generic linear
functional on $\c^N$. Then the number $\mu(f, \pi(\tor))$ of
critical points of $f_{|_{\pi(\tor)}}$ is equal to the degree of
$\pi(\tor)$ as a subvariety of $\c^N$.
\end{Prop}

\begin{proof}
The linear functional $f$ correspond to a Laurant polynomial on
the torus $\tor$ which we denote by $F(x_1, \ldots, x_n)$. A point
$\pi(x), x \in \tor$ is a critical point for $f_{|_{\pi(\tor)}}$
iff $x$ is a critical point of $F$, i.e. if $x$ is a solution of
the system of equations:

$
\begin{cases}
\partial{F} /\partial{x_1} (x) = 0 \\
\cdots \\
\partial{F} /\partial{x_n}(x) = 0 \\
\end{cases}
$

Let $\Delta$ be the Newton polyhedron of $F$ and denote by
$\Delta_i$ the Newton polyhedron of $\partial{F} /\partial{x_i}$.
By Bernstein-Kuchnierenko theorem we know that the number of
solutions of the above system is equal to $ n! \cdot V(\Delta_1,
\ldots, \Delta_n)$, where $V$ denotes the mixed volume of convex
bodies in $\r^n$.

Notice that since the origin is inside the $\Delta$, $\Delta_i$ is
just $\Delta$ shifted in the direction $-e_i$, where $e_1, \ldots,
e_n$ is the standard basis for the Lattice $\z^N$. But mixed
volume is invariant under shifting of the convex bodies so the
number of solutions is equal to $ n! \cdot V(\Delta, \ldots,
\Delta) = n! \cdot Vol(\Delta)$. By Kouchnirenko's theorem the
last number is equal to the degree which proves that the number of
critical points is equal to the degree.
\end{proof}

\subsection{How one can relate the number of critical points and the degree
of an orbit}
As usual suppose a complex Lie group $G$ acts linearly on $\c^N$.
In this section we try to compare the number of critical points of
a generic linear functional $f$ restricted to an orbit $X$ and the
degree of the variety $X$.

Let us recall some basic definitions from the theory of Lie group
actions. For any $x \in \c^N$ there is a linear map $T: \g
\rightarrow T_xO_x$, where $O_x$ denotes the orbit of $x$, defined
by
$$ T(\xi) = \frac{\partial}{\partial{t}}_{|t=0} exp(t \cdot \xi) \cdot x.$$

Let $Vec(\c^N)$ denote the vector space of all vector fields on
$\c^N$. We can get a linear map $v: \g \rightarrow Vec(\c^N), \xi
\mapsto v_\xi$, by defining $v_{\xi}(x) = T(\xi)(x)$. This is the
so-called generating vector field of $\xi$. Each generating vector
field is tangent to the orbits.

As usual let $f$ be a linear functional on $\c^N$. $f$ can be
differentiated along vector fields. Derivative of $f$ along a
vector field $v \in Vec(\c^N)$ at a point $x$ is $df(x)(v)$. Let
$M \subset \c^N$ be a submanifold. If derivative of $f$ at a point
$x \in M$ along any vector field tangent to $M$ is zero then
$df_{M}(x) = 0$, i.e. $x$ is a critical point for $f$ restricted
to $M$.

Now for $\xi \in \g$ consider the hyperplane $$H_{\xi,f} = \{ x
\in \c^N | df(v_\xi (x)) = 0 \}$$ and let $H_f = \bigcap_{\xi \in
\g} H_{\xi,f}$. Let $X$ be an orbit and let $x \in X$. Since
generating vector fields span the tangent spaces to the orbits we
have $df_{|X}(v_{\xi}(x)) = 0, \forall{\xi} \in \g$ iff $x$ is a
critical point of $f_{|X}$. In other words, the critical points of
$f_{|X}$ are exactly the points in $H_{f} \cap X$. Next
proposition tells us what is dimension of $H_{f}$ as a linear
subspace of $\c^N$.

\begin{Prop}
If $d$ is the maximum dimension of the orbits then $dim(H_{f}) = N
- d$, i.e. $H_{f}$ has complementary dimension to the orbits of
maximum dimension.
\end{Prop}
\begin{proof}
Let $x \in \c^N$ be a point such that its orbit $O_x$ has maximal 
dimension equal to $d$. Let $v_1, \ldots, v_d$ be generating 
vector fields of the Lie algebra vectors $\xi_1, \ldots, \xi_d$ 
respectively such that $v_1(x), \ldots, v_d(x)$
give a basis for $T_xO_x$. Since $v_1, \ldots, v_d$ are indepedent
at $x$, they are independent in a neighbourhood $U$ of $x$ and since 
$d$ is the maximum dimension of orbits, we get that for $y \in U$,
$v_1(y), \ldots, v_d(y)$ give a basis for $T_yO_y$. From definition of 
$H_f$ it follows then that $H_f \cap U = H_{\xi_1, f} \cap \cdots H_{\xi_d, f} 
\cap U$.
That is $H_f \cap U$ has codimension $d$, which implies that $H_f$ itself has 
codimension $d$.
\end{proof}

Let us say a few words about the orbits with maximum dimension for
actions of reductive groups. It turns out that for a linear action
of a complex reductive group, almost all the orbits are the same
and one can talk about the notion of \textit{generic orbit}, more
precisely
\begin{Th}[{Vinberg-Popov} p.]
For a linear action of a complex reductive group $G$ on $\c^N$,
there exists a subgroup $S$ of $G$ and an open dense subset $U$ of
$\c^N$, such that stabilizer of any $x \in U$ is conjugate to $S$.
\end{Th}
Such an $S$ is called \textit{stabilizer in general position} and
the orbits in $U$ are called \textit{generic orbits}. It is not
difficult to show that generic orbits are exactly the orbits with
maximum dimension. Hence for a linear action of a reductive group
$dim(H_{f}) = dim(S)$.

Since the degree of a subvariety in $\c^N$ is its number of
intersection points with a generic plane of complementary
dimension, the above discussion then suggests that for an orbit
$X$ with maximum dimension and a generic linear functional $f$,
$deg(X) = \mu(f, X)$. But unfortunately this is not true. The
reason is that while $f$ is generic the plane $H_f$ is not generic
and no matter what $f$ is, $H_f$ and $X$ can always have
intersection points at infinity. What we intend to prove is that
$deg(X) = \mu(f, X)$ in some good cases. The following is a simple
example that $X$ and $H_f$ always have intersection at infinity
and hence $deg(X)$ and $\mu(f,X)$ are different. In this example
the group is not reductive.

\begin{Ex}
Consider the highest weight representation $V_n$ of $\SL$. As
usual, this is realized as the vector space generated by the
symbols $x^n$, $x^{n-1}y$, $\ldots$, $y^n$. Let $B$ be the Borel
subgroup of upper triangular matrices, and $U$ the unipotent
subgroup of triangular matrices with $1$'s on the diagonal. Since
$x^n$ is an eigen vector for the action of $B$, it is the highest
weight vector. Weyl group of $\SL$ has $2$
elements and $w = \left[\begin{matrix} 0 & 1 \\ 1 & 0 \\
\end{matrix} \right]$ is the non-trivial element in it and $y^n =
w \cdot x^n$.

Let us consider the $U$-orbit of $y^n$ which we denote by X
$$ X = U \cdot y^n = \{ (x+uy)^n | \forall u \in \c \}.$$
We identify $V_n$ and $\c^{n+1}$ via the basis $\{x^n, \ldots, y^n
\}$ and represent vectors in $V_n$ by $n+1$-tuples of numbers.
Then
$$ X = U \cdot y^n = \{ (1, u, \ldots, u^n) |\forall u \in \c\}.$$

Obviously $X$ is a subvariety of degree $n$ of $V_n$. Now let us
see what is closure of $B$-orbit of $y^n$ in $\c P^n$. We
represent points in the projective space by homogeneous their
coordinates. Letting $u$ go infinity, it can be easily verified
that the closure of $X$ in $\c P^n$ is
$$ \bar{X} = X \cup \{(0: \ldots : 0 : 1 : 0) \}.$$

Let $f \in V_n^*$ be a linear functional. Take the vector $\xi =
\left[\begin{matrix} 0 & 1 \\ 0 & 0 \\ \end{matrix} \right] \in
\mathcal{b}$, the Lie algebra of $B$. We want to compute the
generating vector field of this vector and the hyperplane $H_{\xi,
f}$ corresponding to $f$ and this vector. From definition we have
$$ v_{\xi}(x) = \frac{\partial}{\partial{t}}_{|t=0} exp(t \cdot \xi) \cdot x.$$

One has $exp(t \cdot \xi) = \left[\begin{matrix} 1 & t \\ 0 & 1 \\
\end{matrix} \right] \in B$ and
\begin{eqnarray*}
exp(t \cdot \xi) \cdot \left[\begin{matrix}x_1\\x_2\\ \vdots \\
x_{n+1} \\ \end{matrix} \right] &=& \left[\begin{matrix} 1 & t &
t^2 & \ldots & t^n \\ 0 & 1 & 2t & \ldots & 2t^{n-1}\\ & \cdots &
\\ 0 & \ldots & 0 & 1 & nt \\ 0 & 0 & \ldots & 0 & 1\\
\end{matrix}
\right] \left[ \begin{matrix} x_1 \\ x_2 \\ \vdots \\ x_{n+1} \\
\end{matrix} \right] \cr
&=& \left[ \begin{matrix} x_1 + t\cdot x_2 + \cdots \\
x_2 + 2t \cdot x_3 + \cdots \\ \vdots \\ x_n + nt \cdot x_{n+1} + \cdots \\ x_{n+1} \\
\end{matrix} \right]. \cr
\end{eqnarray*}
Hence for $x = \left[\begin{matrix}x_1\\x_2\\
\vdots \\ x_{n+1} \\ \end{matrix} \right]$,
$$ v_\xi (x) = \frac{\partial}{\partial{t}}_{|t=0} exp(t \cdot \xi) \cdot \left[\begin{matrix}x_1\\x_2\\
\vdots \\ x_{n+1} \\ \end{matrix} \right] = \left[\begin{matrix}x_2\\2x_3\\ \vdots \\
nx_{n+1}\\0 \\ \end{matrix} \right].$$

Represent $f$ by $(f_1, \ldots, f_{n+1})$ in the dual basis for
$V_n^*$. Then
\begin{eqnarray*}
H_{f, \xi} &=& \{ x \in V_n | f(v_\xi (x)) = 0 \} \cr &=& \{ x \in
V_n | f_1 \cdot x_2 + 2f_2 \cdot x_3 + \cdots + nf_n \cdot x_{n+1}
= 0 \}
\end{eqnarray*}
every element of $X$ is of the form $(1, u, \ldots, u^n)$.
Remembering the formula for $v_\xi$, the critical points
correspond to the solutions of
$$  \sum_{i=1}^{n} i f_i u^i = 0 $$
Notice that the point at infinity $(0: \ldots : 0 : 1 : 0) \in
\bar{X} \setminus X$ satisfies the equation $f_1 \cdot x_2 + 2f_2
\cdot x_3 + \cdots + nf_n \cdot x_{n+1} = 0$ and hence lies on the
closure of $H_{f, \xi}$. Since $deg(X) = n$, $H_{f, \xi}$
intersects $X$ in at most $n-1$ finite points, which are in fact
the critical points of $f_{|X}$.
\end{Ex}

But for torus actions, the number of critical points is equal to
the degree of orbit. As was discussed before, this in fact follows
form the theorems of Bernstein and Kouchnierenko. In case of a
diagonal action of a torus on $\c^N$, any point all whose
components are non-zero has generic orbit.

\begin{Prop}
Let a torus $\tor$ acts on $\c^N$ via a diagonal representation.
Let $X$ be a generic orbit and $f$ a generic linear functional.
Then we have $deg(X) = \mu(f, X)$.
\end{Prop}
\begin{proof}
This was proved during the proof of Bernstien theorem ({}).
\end{proof}

In the next section, we intend to show that this is true in the
more general case of actions with generic spherical orbits, more
precisely
\begin{Th}
Suppose a connected reductive group $G$ acts linearly on $\c^N$
such that generic orbits are closed and spherical. Let $X$ be a
generic orbit. Then for a generic linear functional $f$ on $\c^N$
we have $$ deg(X) = \mu(X, f).$$
\end{Th}

\subsection{A remark on the number of critical points and incidence 
correspondence} 

Let $\overline{X}$ be a smooth subvariety in $\c P^n$ and let $X = 
\overline{X} \cap \c^n$ be its affine part. Let $f$ be a generic 
linear functional on $\c^n$. $deg(\overline{X})$ is, in fact, degree 
of the map $f_{|X}: X \rightarrow \c$. On the other hand, $\mu(X, f)$, the
 number of critical points of $f_{|X}$, can also be realized as the 
 degree of a map but on the so-called \textit{incidence variety}.
 Let ${\c P^n}^*$ be the Grassmanian of $(n-1)$-planes in $\c P^n$.
 
\begin{Def}
The incidence variety $\Phi$ of $\overline{X}$ is
$$ \Phi = \{ (p, H) | p \in \overline{X}, H \in {\c P^n}^* \textup{ and }
T_p \overline{X} \subset H \}. $$
\end{Def}
We have the natural maps $\pi_1: \Phi \rightarrow X, \pi(p, H) = p$ and 
$\pi_2: \Phi \rightarrow {\c P^n}^*, \pi_2(p, H) = H$. $\pi_2(\Phi)$ is 
called the \textit{dual variety} of $\overline{X}$ and can be thought of 
as the variety of all tangent hyperplanes to $\overline{X}$. 
 
\begin{Prop}
$\mu(f, X) = deg(\pi_2)$.
\end{Prop}
\begin{proof}
$x \in \c^n$ is a critical point of $f_{|X} \iff df_{|X}(x) = 0 \iff 
T_x X \subset ker(f)$. Let $ker(f) = H \in {\c P^n}^*$, then $T_x X 
\subset H \iff (x, H) \in \Phi$, i.e. $(x, H) \in {\pi_2}^{-1}(H)$.
\end{proof}

\subsection{The Example of $\SL$}
In this example we look at image of $\SL$ embedded in some space
of matrices and compute its degree, Euler characteristic of
hyperplane sections and the number of critical points of a generic
linear functional on image of $\SL$. We will see that in this
case, the degree is not equal to the number of critical points.

As is known all irreducible representation of $SL_2(\c)$ are $V_n
= Sym(\c^2)^n$ for any natural number $n$. Obviously $dim(V_n) =
n+1$ and as $SL_2(\c)$ is a simple group all its non-trivial
representations are faithful. These irreducible representations
call also be realized as $V_n =$ vector space of all homogeneous
polynomials in $x$ and $y$ of degree $n$, and $SL_2(\c)$ acts by
$g \cdot f(x,y) = f(g^-1(x,y))$. $V_n$ is in fact the highest
weight representation of $\SL$ corresponding to the weight $n$.

We will use the first realization of $V_n$ namely $Sym(\c^2)^n$.
Denote by $\pi$ the homomorphism from $\SL$ to $GL(Sym(\c^2)^n)$.
Let $\{x,y\}$ be a basis for $\c^2$. Then all the monomials $\{x^i
y^j | i,j \geq 0; i+j = n\}$ give a basis for $V_n$ and $\SL$ acts
according to $g \cdot (x^iy^j) = (ax+by)^i(cx+dy)^j)$, where $g =
\left[\begin{matrix}
a & b \\
c & d \\
\end{matrix}
\right] .$

One can compute the matrix of $\pi(g)$ and see that each entry of
this matrix is a homogeneous polynomial of degree $n$ in $a, b, c$
and $d$ and all possible monomials of degree $n$ in $a, b, c$ and
$d$ appear in the entries. Hence the degree of $\pi(\SL)$ will be
the number of solutions of the system of equations

$ \begin{cases}\label{} P_1 (a, b, c, d) = C_1 \\ P_2(a, b, c, d) =
C_2 \\ P_3(a, b, c, d) = C_3 \\ ad-bc = 1
\end{cases}$ 

where each $P_i$ is a generic homogeneous polynomial in $a, b, c$
and $d$.

We will use Kouchnierenko theorem to find the number of solutions
of such a system. Since $ad-bc$ is a non-degenerate quadratic 
form, after a suitable linear change of coordinates (namely
substitute $a$ and $d$ with $a+id$ and $a-id$ respectively and 
substitute $b$ and $c$ with $ib+c$ and $ib-c$ respectively), it can be 
put in the form $a^2 + b^2 + c^2 + d^2$. It is evident that after
this change of coordinates the polynomials $P_i$ will remain a 
generic polynomial of degree $n$. With abuse of notation we 
denote the polynomials $P_i$ after this change of coordinates 
again by $P_i$. Let $\Delta$ denote the 
polytope in $\r^4$ whose vertices are the origin and the standard 
basis vectors. $4$ dimensional volume of $\Delta$is $1/4!$. Then the Newton polyhedron of $a^2 + b^2 + c^2 + d^2$
is $2\Delta$ and the Newton polyhedron of $P_i(a, b, c, d) - C_i$ 
is equal to $n\Delta$. From the Bernstein-Kouchnierenko theorem 
(Theorem ~{}) it follows that the number of the solutions of 
the above system is equal to
\begin{eqnarray*}
4! V(2\Delta, n\Delta, n\Delta, n\Delta) &=&
4! \. 2n^3 \. Vol(\Delta) \cr
&=& 4! \. 2n^3 \. 1/4! \cr
&=& 2n^3
\end{eqnarray*}


That is, the degree of $\pi(\SL)$ in $M(n+1, \c)$, the vector
space of $(n+1) \times (n+1)$ matrices, is equal to $2n^3$.

On the other hand, Euler characteristic, denoted by $\chi$, of a
hyperplane section of $\pi(\SL)$ is equal to Euler characteristic
of $\{(a,b,c,d)|ad-bc=1; P(a,b,c,d)=C\}$, where $P$ is a
generic homogeneous polynomial of degree $n$ and $C$ is some
generic constant number.As before with a linear change of 
coordinates we can work with $a^2 + b^2 + c^2 + d^2 = 1$ instead of 
$ad-bc = 1$. With abuse of notation, we denote the polynomial $P$
after the change of coordinate by $P$ again. So it sufices to
find $\chi(Z)$ where
$$ Z = \{(a,b,c,d)|a^2+b^2+c^2+d^2=1; P(a,b,c,d)=C\}.$$

Recall that Euler characteristic is an additive function i.e. if
$X = A \cup B$ and $A \cap B = \emptyset$ then $\chi(X) = \chi(A)
= \chi(B)$. We partition Z = $\{(a,b,c,d)|ad-bc=1; P(a,b,c,d)=C\}$ 
into pieces whose Euler characteristic can be computed using 
Bernstein's theorem.

Let $S \subset \{ a,b,c,d \}$ and define $Z_S \subset Z$ to be the 
set of points in $Z$ that exactly their coordinates belonging to 
$S$ are zero, e.g. $Z_\emptyset =  
\{(a,b,c,d)|a^2+b^2+c^2+d^2=1; P(a,b,c,d); a,b,c,d \neq 0 \}$ and
$Z_a = \{(a,b,c,d)|a^2+b^2+c^2+d^2=1; P(a,b,c,d); a,b,c,d \neq 0 \}$.
Then $Z$ is the disjoint union of $Z_S$'s $(S \subset \{a,b,c,d\}$.
From additivity of the Euler characteristic we obtain
$$\chi(Z) = \sum_{S \subset \{a,b,c,d\}} \chi(Z_S).$$
It is evident that 
$$\chi(Z_a) =\chi(Z_b) =\chi(Z_c) = \chi(Z_d).$$
and 
$$\chi(Z_{ab}) =\chi(Z_{bc}) =\chi(Z_{cd}) = \chi(Z_{ad}) = \chi(Z_{ac})
 = \chi(Z_{bd}).$$
and 
$$\chi(Z_{abc}) =\chi(Z_{bcd}) =\chi(Z_{acd}) = \chi(Z_{abd}).$$
Hence
$$\chi(Z) = \chi(Z_\emptyset) + 4 \chi(Z_a) + 6 \chi(Z_{ab}) + 
4 \chi(Z_{abc}).$$
Next we compute each of the terms separately using generalized version of 
Bernstein's theorem. In the following  
multiplicative notation refers to mixed volume of polyhedra.
\begin{enumerate}
\item $\chi(Z_\emptyset)$ is the $4$th degree term in the expansion of 
$$4!(2\Delta - (2\Delta^2) + (2\Delta)^3)(n\Delta-(n\Delta)^2+(n\Delta)^3)$$
That is 
\begin{eqnarray*}
4!(2n^3\Delta^4+4n^2\Delta^4+8n\Delta^4) 
&=& 4! (2n^3+4n^2+8n) Vol_4(\Delta) \cr
&=& 2n^3+4n^2+8n
\end{eqnarray*}
\item Let $\Delta'$ be the polyhedron in $\r^3$ whose vertices are 
the origin and the standard basis vectors. $3$ dimensional volume of 
$\Delta'$ is $1/3!$. $\chi(Z_a)$ is then the $3$rd degree term in the 
expansion of 
$$3!(2{\Delta'} - (2{\Delta'}^{2}))(n{\Delta'}-(n{\Delta'})^{2})$$
That is 
\begin{eqnarray*}
3!(-2n^2{\Delta'}^{3}-4n{\Delta'}^{3}) 
&=& 3! (-2n^2-4n) Vol_3({\Delta'}) \cr
&=& -2n^2-4n
\end{eqnarray*}
\item Let $\Delta"$ be the polyhedron in $\r^2$ whose vertices are 
the origin and the standard basis vectors. $2$ dimensional volume of 
$\Delta'$ is $1/2!$. $\chi(Z_ab)$ is then the $2$nd degree term in the 
expansion of 
$$2! \. 2\Delta" \. n\Delta"$$
That is 
\begin{eqnarray*}
2! \. 2n {\Delta"}^2) 
&=& 2! \. 2n Vol_2(\Delta") \cr
&=& 2n
\end{eqnarray*}

\item Finally note that for generic polynomial $P$, $Z_{abc}$ is 
empty since $a^2+b^2+c^2+d^2=1$ and $a=b=c=0$ implies $d=\pm 1$ and 
$\pm 1$ is not a root of a generic polynomial $P(0,0,0,d) - C$.
\end{enumerate}

Putting every thing together, we obtain
\begin{eqnarray*}
\chi(Z) &=& \chi(Z_\emptyset) + 4 \chi(Z_a) + 6 \chi(Z_{ab}) + 
4 \chi(Z_{abc}) \cr
&=& (2n^3+4n+8n) + 4(-2n^2-4n) + 6(2n) + 4\.0 \cr
&=& 2n^3 - 4n^2 +4n
\end{eqnarray*}

Hence Euler characteristic of a generic hyperplane section of
$\pi(\SL)$ in $M(n+1, \c)$, the vector space of $(n+1) \times
(n+1)$ matrices, is equal to $2n^3-4n^2+4n$. Of course this is 
not in general equal to $2n^3$, although for $n=1$ 
,i.e. the natural representaiton of $\SL$ the two numbers are 
coincide
$$ 2.1^3 - 4 \. 1^2 + 4 \. 1 = 2 \. 1^3 = 2. $$
More generally, we will show that, this is the case for the natural
representation of $SL(n, \c)$. In general, as is the case in this example,
the main terms of the formulae for degree and the number of critical points
are the same. We will see this in the last section (section ~{}).

\section{Actions with spherical orbits and the main theorem}

In this section we intend to prove that for a group action with
generic spherical orbits, the number of critical points of a
generic functional on a generic orbit is equal to the degree of
the orbit. Before we discuss the proof we need some preliminaries
regarding spherical varieties.

\subsection{Some preliminaries form the theory of spherical varieties}

Recall that a $G$-variety is a spherical variety if a Borel
subgroup of $G$ has a dense open orbit.

An important property of the spherical varieties is that one can
approach any point in the closure of an orbit with a one parameter
subgroup. Next theorem is even a stronger result. The proof can be
found in {Brion-Luna-Vust}. First we need the definition of an
equivariant embedding.

\begin{Def}
Let $G/H$ be a homogeneous space. An \textit{equivariant
embedding} of the homogeneous space $G/H$ is a pair $(X, i)$,
where $X$ is a variety with a $G$ action and $i$ is a
$G$-equivarint embedding $i: G/H \rightarrow X$ with $i(G/H)$ an
open dense subset of $X$.
\end{Def}

\begin{Th}[Brion-Luna-Vust]
Let $G/H$ be a spherical homogeneous space. Let $x \in G/H$ then
there exists a torus $T_x \subset G$, not necessarily unique, such
that in any equivariant embedding $Y$ of $G/H$, closure of the
$T_x$-orbit of $x$ intersects all the $G$-orbits.
\end{Th}

Let us see how the above theorem implies that one can approach any
point in the closure of any orbit with a one parameter subgroup:
suppose $y \in \bar{G/H}$. By the above theorem, let $g \cdot y$
belong to the intersection of closure of the $T_x$-orbit of $x$
and $G$-orbit of $y$. Hence there is a one dimensional subtorus
$S$ (one parameter subgroup) of $T_x$ such that $S \cdot x$
approaches $g \cdot y$. Now if $S' = g \cdot S \cdot g^{-1}$ and
$x' = g \cdot x$ then $S'$ is a one parameter subgroup of $G$ and
$S' \cdot x'$ approaches $y$.

\begin{Cor}
Any spherical variety consists of a finite number of $G$-orbits.
\end{Cor}
\begin{proof}
Closure of $T_x$-orbit of $x$, is a toric variety. It is
well-known that a toric variety consists of a finite number of
torus orbits. But closure of this $T_x$-orbit intersects all the
$G$-orbits and hence there should exist a finite number of
$G$-orbits.
\end{proof}

In fact even the number of Borel orbits is finite

\begin{Th}
In a spherical variety, for a Borel subgroup $B$, there are only a
finite number of $B$ orbits.
\end{Th}

For the proof, look at {Brion-Luna-Vust}.

\subsection{The main thoerem}
Now we are in the position to state and prove the theorem
regarding the degree of a generic spherical orbit and the number
of critical points.

\begin{Th}
Let $G$ acts linearly on $\c^N$ such that generic orbits are
spherical and closed. Let $X$ be a generic orbit and let $f$ be a
generic linear functional on $\c^N$. Then
$$ deg(X) = \mu(X, f).$$
\end{Th}
\begin{proof}

As we saw in section ~{}, $p \in X$ is a critical point of $f_{|X}$
iff for any generating vector field $v_{\xi} (\xi \in \g)$ on $X$, 
$f(v_{\xi}(x)) = 0$, i.e. $p$ belongs to the hyperplane 
$H_f = \{ x \in \c^n | f(v_{\xi}(x) = 0 \}$.
By Proposition ~{}, $H_f$ has codimension equal to the dimension of 
$X$. Hence intersection number of $\overline{X}$ and $\overline{H}$
(that is closures of $X$ and $H$ in $\c P^n$) is equal to  $deg(X)$.
If $f$ is generic, the critical points of $f_{|X}$ are non-degenerate
and hence $\mu(f, X) = \sharp (H_f \cap X)$. So we only need to show that 
$\overline{X}$ and $\overline{H_f}$ do not have any intersection at infinity.,
i.e. $X \cap H_f = \overline{X} \cap \overline{H_f}$. 

We proceed the proof by contradiction. Suppose there exists $z \in  
\overline{X} \cap \overline{H_f} \setminus X \cap H_f$. Let $B$ be the 
Borel subgroup of $G$ with a dense orbit. Choose coordinates in $\c^n$ 
such that $B$ acts by upper triangular matrices. Define 
$$ V_i = \{ (x_1 : x_2 : \ldots : x_i : 1 : 0 : \ldots : 0 ) \in \c P^n \}. $$
where we have used homogenuous coordinates to represent points in the 
projective space. Then $V_i \cong \c^i$ and we have the cell 
decomposition $\c P^n = V_n \cup \cdots \cup V_1$. 
  
Let the linear functional $f$ be given by $f(x) = \sum_{i=1}^{n} f_ix_i$, 
where $x = (x_1, \ldots x_n) \in \c^n$. Define a linear functional on $V_i$
by 
$$ \tilde{f}_i (x) = \sum_{j=1}^{i} f_jx_j$$ 
where $x = (x_1 : x_2 : \ldots : x_i : 1 : 0 : \ldots : 0 ) \in V_i$.  
Suppose $z = (z_1 : z_2 : \ldots : z_k : 1 : 0 : \ldots : 0 ) \in V_k$.
We intend to show:

\textit{ $\sum_{j=1}^{k+1} f_jz_j = 0$ (that is $\tilde{f}_k(z) = -f_{k+1}$) 
and $z$ is a critical point $\tilde{f}_k$ restricted to $B \. z$, the 
Borel orbit of $z$ (that is $-f_{k+1}$ is a critical value).}

But by Sard's theorem, for generic $f$, $-f_{k+1}$ is not a critical 
value. Since $\overline{X}$ consists of a finite number of Borel orbits 
$-f_{k+1}$ is not a critical value for $\tilde{f}_i$ restricted to a Borel 
orbit. Hence there is no point 
$z \in \overline{X} \cap \overline{H_f} \setminus X \cap H_f$.

It only remains to prove the above claim. Next lemma proves the 
first claim, i.e. if $z \in \overline{X} \cap \overline{H_f}$ then 
$\sum_{j=1}^{n} f_jz_j = \tilde{f}_k(z) + f_{k+1} = 0$. 

Notice that from the theorem ~{}, we can approach $z$ with a 
one-parameter subgroup $\alpha$ of $G$.

Let $GL(n, \c)$ acts in the usual way on $\c^n$.
Let $z \in \c P^n \setminus \c^n$ be a point at infinity. 
Suppose there is a point 
$x \in \c^n$ and a one-parameter subgroup $exp(t \. \alpha), \alpha 
\in \mathfrak{gl}(n, \c)$, with $\alpha$ diagonalizable, such that 
$exp(t \. \alpha) \. x$ converges to $z$ as $t$ goes to infinity. 
Denote the generating vector field of $\alpha$ on $\c^n$ by $v_{\alpha}$. 
Let $f$ be a linear functional on $\c^n$ and as before $H_{f, \alpha}$ be 
the hyperplane $\{y | df(v_{\alpha}(y) = 0 \}$, we have

\begin{Lem}
$z \in \overline{H_{f, \alpha}}$ implies that 
$z \in \{ y | f(y) = 0 \}$. Roughly speaking: if we can approach a 
point $z$ at infinity with a one-parameter subgroup $\alpha$ and if 
as we approach $z$ the derivative of $f$ along $\alpha$ is zero, 
then $f$ lies on the hyperplane at infinity defined by $f$. 
\end{Lem}

\begin{proof}
The proof is based on direct calculation of $v_{\alpha}$.
Choose coordinates in $\c^n$ such that $\alpha$ becomes diagonal. So let us 
assume that $\alpha = (\alpha_1, \ldots, \alpha_n)$. For $y = (Y_1, 
\ldots, y_n)$ we have
\begin{eqnarray*}
v_{\alpha} (y) &=& \frac{\partial}{\partial s}_{|s=0} exp(s \alpha)\. y  
\cr 
&=&  \frac{\partial}{\partial s}_{|s=0} (e^{\alpha_1}\.s, \ldots , 
e^{\alpha_n}\.s ) \. (y_1, \ldots, y_n) \cr
&=& (\alpha_1 y_1, \ldots, \alpha_n y_n)
\end{eqnarray*}

Now let $x = (x_1, \ldots, x_n)$ then 
\begin{eqnarray*}
exp(t \alpha) \. x &=& (e^{t \. \alpha_1} \. x_1, \ldots, 
e^{t \. \alpha_n} \. x_n) \cr
&=& (e^{t \. \alpha_1} \. x_1 : \ldots : 
e^{t \. \alpha_n} \. x_n : 1) \cr
&=& (e^{t \. (\alpha_1 - m)} \. x_1 : \ldots : 
e^{t \. (\alpha_n - m)} \. x_n : e^{-t \. m}) \cr
\end{eqnarray*}

where $m = min(\alpha_1, \ldots, \alpha_n)$. Note that all 
$\alpha_i - m \geq 0$ and hence $e^{t \. (\alpha_i - m)}$ is either 
$1$ or approaches $0$ as $t \rightarrow \infty$. Also since 
$\lim_{t \rightarrow \infty} exp(t \alpha) \. x = z \notin \c^n$
then $\lim_{t \rightarrow \infty} e^{-t\. m} = 0$ so $m > 0$ and we have
\begin{displaymath}
z_i = \left \{ \begin{array}{ll} 0 & \textrm{if $\alpha_i \neq m$} \\
x_i & \textrm{if $\alpha_i = m$} \end{array} \right.
\end{displaymath}
Let $f(y) = \sum_{i=1}^{n}f_iy_i$. Suppopse 
$z \in overline{H_{f, \alpha}}$,
then $\sum_{i=1}^{n}f_i \alpha_i z_i = 0$ where 
$z = (z_1 : \ldots : z_n : 0)$. But $z_i = 0$ if $\alpha_i \neq m$, hence
\begin{eqnarray*}
0 &=& \sum_{i=1}^{n}f_i \alpha_i z_i \cr
&=& \sum_{\{i | \alpha_i = m \}} f_i \alpha_i z_i \cr
&=& m \. \sum_{\{i | \alpha_i = m \}} f_i z_i \cr
&=& m \. \sum_{i=1}^{n}f_i z_i \cr
\end{eqnarray*}

Since $m > 0$, this implies that $\sum_{i=1}^{n}f_i z_i = 0$
as was required.
\end{proof}

Next we will show that $z \overline{X} \cap \overline{H_f} \setminus 
X \cap H_f$ implies that $z$ is a critical point of 
$\tilde{f}_k{_|B \. z}$. For this, we neeed to calculate the generating vector 
fields of Borel elements on $V_i$'s

Let $B \subset GL(n, \c)$ be the subgroup of all upper triangular matrices.
Natural action of $GL(n, \c)$ on $\c^n$ extends to $\c P^n$ by
letting it act trivially on the last homogenous component. One can 
easily see that

\begin{Prop}
The action of $B$ on $\c P^n$ respects the cell decomposition 
$\c P^n = V_n \cup \cdots \cup V_1$. 
\end{Prop}

Let $\xi \in \b$, 
the Lie algebra of $B$ and $exp(t \. \xi) = [b_{ij}(t)]$ the corresponding 
one-parameter subgroup. The diagonal elements of $exp(t \. \xi)$ are 
homomorphisms $b_{ii}(t) = e^{c_i t}$. Next simple lemma gives a formula for 
the generating vector field of $\xi$ on the cells $V_i$.

\begin{Lem} 
Let 
$z \in V_k = \{ (x_1:x_2:\ldots:x_k:1:0:\ldots:0) \in \c P^n \}$ and let 
$v_{\xi}$ denote the generating vector field of $\xi$ on $V_k$ then 

$$v_{\xi}(z) = (\sum_{j=1}^{k} b'_{1j}(0)z_j - c_{k+1}z_1 :  
\sum_{j=2}^{k} b'_{2j}(0)z_j - c_{k+1}z_2 \\ :\ldots, b'_{kk}(0)z_k - 
c_{k+1}z_k : 1 : 0 : \ldots : 0) $$

In particular if $k = n$, i.e. generating vector field on $\c^n$ 

$$v_{\xi}(z) = (\sum_{j=1}^{n} b'_{1j}(0)z_j,  
\sum_{j=2}^{n} b'_{2j}(0)z_j ,\ldots, b'_{nn}(0)z_k) $$
\end{Lem}
\begin{proof}
Just from the definition we have
$$ v_{\xi}(z) = \frac{\partial}{\partial{t}}_{|t=0} exp(t\xi) \. z.$$  
Also from the definition of the action we have

\begin{eqnarray*} 
exp(t\xi) \. z &=& (\sum_{j=1}^{k} b_{1j}(t)z_j:  
\sum_{j=2}^{k} b_{2j}(t)z_j :\ldots, b_{{k+1}{k+1}}(t): 0 : 
\ldots : 0). \cr
&=& (\sum_{j=2}^{k} \frac{b_{2j}(t)z_j}{b_{{k+1}{k+1}}(t)} :\ldots : 1 : 0 : 
\ldots : 0). \cr
\end{eqnarray*}
Differentiating with respect to $t$ we get
$$ \frac{\partial}{\partial{t}}_{|t=0} exp(t\xi) \. z =
(\sum_{j=1}^{k} \partial{(b_{1j}'(0)\. b_{{k+1}{k+1}}(0) - b_{1j}(0) 
\. b_{{k+1}{k+1}}'(0))}{b_{{k+1}{k+1}}(0)^2} : \\ 
\ldots : 1 : 0 : \ldots : 0).$$
but $b_{{k+1}{k+1}}(0)=1, b_{ij}(0) = \delta_{ij}$ and $b'_{{k+1}{k+1}}(0)
= c_{k+1}$ so 
$$v_{\xi}(z) = (\sum_{j=1}^{k} b'_{1j}(0)z_j - c_{k+1}z_1 :  
\sum_{j=2}^{k} b'_{2j}(0)z_j - c_{k+1}z_2 \\ :\ldots, b'_{kk}(0)z_k - 
c_{k+1}z_k : 1 : 0 : \ldots : 0) $$
\end{proof}

Having every thing in hand, we know see that $z$ is a critical 
point of $\tilde{f}_{k_{|B\.z}}$.
 
Let $v_{\xi}$ be the generating vector field $\xi \in \b$ on $\c^n$.
From the above lemma (lemma ~{})
$$v_{\xi}(y) = (\sum_{j=1}^{n} b'_{1j}(0)y_j,  
\sum_{j=2}^{n} b'_{2j}(0)y_j ,\ldots, b'_{nn}(0)y_k) $$
where $y = (y_1, \ldots, y_n)$.
Since $\forall \xi, z \in 
\overline{H_{\xi, f}} = \overline{ \{ y|f(v_{\xi}(y)=0 \} }$ we have

\label{ghaf} 
$$\sum_{i=1}^{n}\sum_{j=i}^{n}f_i\.b'_{ij}(0)\.z_j = 0$$.

On the other hand, since we can approach $z$ with a 
one parameter subgroup from lemma ~{} it follows that 
$z \in \overline{\{ y|f(y) = 0 \}}$, i.e. 
$\sum_{i=1}^{n} f_iz_i = 0$.
Subtract $c_{k+1}\. \sum_{i=1}^{n}f_iz_i$ from the left side
of ~{ghaf} we get
\begin{eqnarray*}
0 &=& \sum_{i=1}^{n}\sum_{j=i}^{n}f_i (b'_{ij}(0)\.z_j - c_{k+1}z_i). \cr
&=& \sum_{i=1}^{k+1} f_i \sum_{j=i}^{k+1} 
(b'_{ij}(0)\.z_j - c_{k+1}z_i). \cr
&=& \sum_{i=1}^{k} f_i \sum_{j=i}^{k} (b'_{ij}(0)\.z_j - c_{k+1}z_i) 
+ f_{k+1} \. (b'_{{k+1}{k+1}}(0)z_{k+1} - c_{k+1}x_{k+1}). \cr
&=& \sum_{i=1}^{k} f_i \sum_{j=i}^{k}(b'_{ij}(0)\.z_j - c_{k+1}z_i). \cr
&=& v_{\xi}(z)
\end{eqnarray*}
Since $v_{\xi}(z) = 0$ for all $\xi \in \b$ we conclude that $z$
is a critical point of $\tilde{f}_k{_|B \. z}$. This finishes the 
proof of the main theorem.
\end{proof}

Finally, from theorem {} we then obtain

\begin{Cor}[the main theorem]
Let $G$ acts linearly on $\c^N$ such that generic orbits are
spherical and closed. Let $X$ be a generic orbit and let $f$ be a
generic linear functional on $\c^N$ and $c$ a generic complex
number. Then
$$ \chi(f^{-1}(c) \cap X) = \chi(X) + (-1)^{dim(X)+1} \cdot
deg(X).$$
\end{Cor}

\begin{Ex}[the main example]
Let $G = \SL$ consider its irreducible representations 
$V_n = Sym^2(\c^n)$. This gives an embedding $\pi: \SL \hookrightarrow 
M(n+1, \c)$.
As we mentioned before, $G \times G$ acts on $M(n+1, \c)$ by left-right
multiplication and $\pi(G)$ is a spherical orbit. Unfortunately, one 
can show that it is a generic orbit of this action only for $n=2$ which 
correspond to natural representation of $\SL$. According to our
previous calculation of the number of critical points and the degree
for $\SL$ embedded in $M(n+1, \c)$ in section ~{SL2}, $degree$ is
equal to $2n^3$ while the number of critical points $\mu(f, \SL) = 
4n^3-6n^2+4n$. In case of $n=1$ these two numbers are equal: $4-6+4 = 2 $. 
Which agrees with the main theorem 
~{main-theorem}, since $\SL$ is a closed generic spherical orbit of
the action of $\SL \times \SL$ on $M(2, \c)$ (it is closed simply because
it is given by $det = 1$).

More generally one can see that for $G = SL(n, \c)$, action of 
$G \times G$ on $M(n, \c)$ is an action with closed generic spherical 
orbits and $SL(n, \c)$ itself is a generic spherical orbit. We can  
verify this as follows: let $a \in M(n, \c)$ be an invertible element.
It is obvious that the $G \times G$ orbit of $a$ is $\{ x | det(x) = det(a) \}$.
which is a hypersurface in $M(n, \c)$. The $G \times G$-Stabilizer of 
$a$ is also $\{ (g, m\.g\.m^{-1} | g \in G \}$. Now let $b$ be another 
invertible matrix. One can easily see that $x\. Stab(a) \.x^{-1} = Stab(b)$ 
where $x = (1, ba^{-1}) \in G\times G$. So all the orbits of invertible 
elements are generic orbtits and of course $SL(n, \c)$ is the orbit of identity.

Noting that $dim(SL(n, \c) = n^2-1$ and $\chi(SL(n, \c) = 0$ (cf. thereom 
~{}), the main theorem in this case then 
implies that
\begin{Th}
The Euler characteristic of a generic hyperplane section of 
$SL(n, \c) \subset M(n, \c)$ is equal to ${-1}^{n^2} \. n$.
\end{Th}
\end{Ex}

\subsection{Classification of modules with spherical orbits}

In {Arzhantsev} has classified all the representations with
generic spherical orbits. We quote the main results of his paper
here in this subsection.

Let $G$ be a connected reductive algebraic group over an
algebraically closed field $K$ of zero characteristic, and $G^s$
denote the maximal connected semisimple subgroup of $G$.

\begin{Def}
Let $X$ be an irreducible algebraic variety. We shall say that an
action $G:X$ is {\it an action with generic spherical orbits} if
there exists an open subset $X_0\subset X$ such that for any $x\in
X_0$ the orbit $Gx$ is spherical.
\end{Def}

In ~[Corollary~1]{Arzhantsev1} it is shown that if a
$G$-module $V$, is a module with generic spherical orbits, then in
fact all $G$-orbits are spherical.

Below we list some basic facts, without proof, about actions with
spherical orbits.

(1)\ Any trivial $G$-action is an action with spherical orbits.

(2)\ Suppose that for an action $G:X$ a stabilizer in general
position exists, see~[sec.~7.3]{Vinberg-Popov}. (This is
always the case for linear actions.) Denote this subgroup by $H$.
The action $G:X$ is an action with spherical orbits iff $H$ is a
spherical subgroup of $G$.

(3)\ Rosenlicht's theorem ~{Rosenlicht}, implies that an
action $G:X$ is an action with spherical orbits iff
$K(X)^G=K(X)^B$, where $K(X)^G$ and $K(X)^B$ denote the subfields
of rational function on $X$ invariant under $G$ and $B$
respectively.

(4)\ A $G$-module $V$ is called a \textit{spherical module} if it
is spherical as a $G$-variety, i.e. if a Borel subgroup has a
dense open orbit in $V$. It is shown that any module with
spherical orbits can be realized as a spherical module after an
extension of the group $G$ by a central torus.

(5)\ Let $G_1:X_1$ and $G_2:X_2$ be actions with spherical orbits.
Then the action $(G_1\times G_2)\,:\,(X_1\times X_2)$ is an action
with spherical orbits.

Because of the above fact (5), to classify the modules with
spherical orbits, it suffices to consider only
\textit{indecomposable modules}

\begin{Def}
A $G$-module $V$ is {\it indecomposable} if there exists no proper
decompositions $G^s=G^s_1\times G^s_2$ and $V=V_1\oplus V_2$ such
that $(g_1, g_2)(v_1, v_2)=(g_1v_1, g_2v_2)$ for any $g=(g_1,
g_2)\in G^s$ and any $v=(v_1, v_2)\in V$.
\end{Def}

Yet, there is another way one can obtain modules with spherical
orbits without any cost, i.e. extension of the group $G$ by a
torus

\begin{Def}
We say that a $G'$-module $V$ is obtained from a $G$-module $V$
{\it by a torus extension} if there exists a torus $T$ acting on
$V$ such that $T$- and $G$-actions commute and $G'=TG$.
\end{Def}

It is clear that any $G$-module $V$ is obtained by a torus
extension from the $G^s$-module $V$.

\begin{Lem} \label{lll}
Suppose that $V$ is a $G$-module with spherical orbits and a
$G'$-module $V$ is obtained from this module by a torus extension.
Then $V$ is a $G'$-module with spherical orbits.
\end{Lem}

\begin{proof}
Let $H$ be the generic isotropy subgroup for the action $G:V$. By
assumption, $H$ is spherical in $G$. Then any subgroup of $G'$
containing $H$ is spherical in $G'$. Hence a generic isotropy
subgroup for the $G'$-module $V$ is spherical.
\end{proof}

Now we are ready to state the classification theorem for
$G$-modules with spherical orbits

\begin{Th} \label{mn}
All indecomposable $G$-modules with spherical orbits are either
indicated in Tables 1-3 or are obtained from the indicated modules
by a torus extension.
\end{Th}

\newpage

\centerline{\bf Table 1}

\scriptsize
$$
\begin{array}{||c|c|c|c|c|c||}
\hline
 & & & & & \\
 & G & \mbox{weights} & \dim V & {\mathcal H} &
 \mbox{codim} \\
 & & & & & \\
\hline
 & & & & & \\
 \ 0 \ & \ \{ e\} \ & \ 0 & 1 \ & \ 0 \ & \ 1 \ \\
 & & & & & \\
\hline
 & & & & & \\
 \ 1 \ & \ SL(n) \ & \ \phi_1 & n \ & \ A_{n-2}+R_{n-1} \ & \ 0 \ \\
 & & & & & \\
\hline
 & & & & & \\
 \ 2 \ & \ \Lambda^2 SL(2n) \ & \ \phi_2 & 2n^2-n \ & \ C_n \ & \ 1 \ \\
 & & & & & \\
\hline
 & & & & & \\
 \ 3 \ & \ \Lambda^2 SL(2n+1) \ & \ \phi_2 \ & \ 2n^2+n \
 & \ C_n+R_{2n} \ & \ 0 \ \\
 & & & & & \\
\hline
 & & & & & \\
 \ 4 \ & \ S^2 SL(2n) \ & \ 2\phi_1 \ & \ 2n^2+n \ & \ D_n \ & \ 1 \ \\
 & & & & & \\
\hline
 & & & & & \\
 \ 5 \ & \ \ S^2 SL(2n+1) \ \ & \ \ 2\phi_1 \ \ & \ \ 2n^2+3n+1 \ \ & \ B_n \
 & \ 1 \ \\
 & & & & & \\
\hline
 & & & & & \\
 \ 6 \ & \ SO(2n) \ & \ \phi_1 \ & \ 2n \ & \ B_{n-1} \ & \ 1 \ \\
 & & & & & \\
\hline
 & & & & & \\
 7 & SO(2n+1) & \phi_1 & 2n+1 & D_n & 1 \\
 & & & & & \\
\hline
 & & & & & \\
 8 & Spin(7) & \phi_3 & 8 & G_2 & 1 \\
 & & & & & \\
\hline
 & & & & & \\
 9 & Spin(9) & \phi_4 & 16 & B_3 & 1 \\
 & & & & & \\
\hline
 & & & & & \\
 10 & Spin(10) & \phi_4 & 16 & B_3+R_8 & 0 \\
 & & & & & \\
\hline
 & & & & & \\
 \ 11 \ & Sp(2n) & \phi_1 & 2n & \ \ C_{n-1}+R_{2n-1} \ \ & 0 \\
 & & & & & \\
\hline
 & & & & & \\
 12 & G_2 & \phi_1 & 7 & A_2 & 1 \\
 & & & & & \\
\hline
 & & & & & \\
 13 & E_6 & \phi_1 & 27 & F_4 & 1 \\
 & & & & & \\
\hline
\end{array}
$$

\newpage

\normalsize
\centerline{\bf Table 2 }
\tiny

$$
\begin{array}{||c|c|c|c|c|c||}
\hline
 & & & & & \\
 & G & \mbox{weights} & \dim V & {\mathcal H} &
 \mbox{codim} \\
 & & & & & \\
\hline
 & & & & & \\
 14 & SL(2)\times K^* &
 \phi_1\otimes\epsilon+\phi_1\otimes\epsilon^{-1} & 4 & t_1 & 1 \\
 & & & & & \\
\hline
 & & & & & \\
 15 & SL(n)\times K^*, \ n>2 &
 \phi_1\otimes\epsilon^a+\phi_1\otimes\epsilon^b, \ a\ne b
 & 2n & A_{n-3}+t_1+R_{2(n-2)} & 0 \\
 & & & & & \\
\hline
 & & & & & \\
 16 & SL(n), \ n>2 &
 \phi_1+\phi_{n-1} & 2n & A_{n-2} & 1 \\
 & & & & & \\
\hline
 & & & & & \\
 17 & SL(2n+1) & \phi_1+\phi_2 & (2n+1)(n+1) & C_n & 1 \\
 & & & & & \\
\hline
 & & & & & \\
 18 & SL(2n+1)\times K^* &
 \phi_1\otimes\epsilon^a+\phi_{2n-1}\otimes\epsilon^b, a\ne nb  &
 (2n+1)(n+1) & C_{n-1}+t_1+R_{2(2n-1)} & 0 \\
 & & & & & \\
\hline
 & & & & & \\
 19 & SL(2n) & \phi_1+\phi_2 & n(2n+1) & C_{n-1}+R_{2n-1} & 1 \\
 & & & & & \\
 & & \phi_1+\phi_{2n-2} & & & \\
 & & & & & \\
\hline
 & & & & & \\
 20 & SO(8) & \phi_1+\phi_3 & 16 & G_2 & 2 \\
 & & & & & \\
\hline
 & & & & & \\
 21 & Sp(2n)\times K^* &
 \phi_1\otimes\epsilon+\phi_1\otimes\epsilon^{-1} &
 4n & C_{n-1}+t_1 & 1 \\
 & & & & & \\
\hline
 & & & & & \\
 22 & SL(n)\times SL(m), \ n>m & \phi_1\otimes \phi_1 & nm &
 A_{n-m-1}+A_{m-1}+R_{nm-m^2} & 0 \\
 & & & & & \\
\hline
 & & & & & \\
 23 & SL(n)\times SL(n) & \phi_1\otimes \phi_1 & n^2 & A_{n-1} & 1 \\
 & & & & & \\
\hline
 & & & & & \\
 24 & SL(2)\times Sp(2n) & \phi_1\otimes \phi_1 &
 4n & C_{n-1}+A_1 & 1 \\
 & & & & & \\
\hline
 & & & & & \\
 25 & SL(3)\times Sp(2n)\times K^*, \ n>1
 & \phi_1\otimes\phi_1\otimes\epsilon &
 6n & C_{n-2}+A_1+t_1+R_{2n-1} & 0 \\
 & & & & & \\
\hline
 & & & & & \\
 26 & SL(4)\times Sp(4) & \phi_1\otimes \phi_1 &
 16 & C_2 & 1 \\
 & & & & & \\
\hline
 & & & & & \\
 27 & SL(n)\times Sp(4), \ n>4 & \phi_1\otimes \phi_1 &
 4n & A_{n-5}+C_2+R_{4(n-4)} & 0 \\
 & & & & & \\
\hline
\end{array}
$$

\newpage

\normalsize
\centerline{\bf Table 3 } 
\scriptsize

$$
\begin{array}{||c|c|c|c|c||}
\hline
 & & & & \\
 & G & \mbox{weights} & \dim V &
 \mbox{codim} \\
 & & & & \\
\hline
 & & & & \\
 28 & SL(n)\times SL(n)\times K^* &
 \phi_1\otimes\epsilon+ \phi_1\otimes\psi_1 &
 n(n+1) & 1 \\
 & & & & \\
  & & \phi_1\otimes\epsilon+\phi_{n-1}\otimes\psi_{n-1} & & \\
 & & & & \\
\hline
 & & & & \\
 29 & SL(n+1)\times SL(n)\times K^* &
 \phi_1\otimes\epsilon^n+\phi_1\otimes\psi_1\otimes\epsilon^{-1} &
 (n+1)^2 & 1 \\
 & & & & \\
\hline
 & & & & \\
 30 & SL(n+1)\times SL(n)\times K^*\times K^*, \ n>1 &
 \phi_1\otimes\epsilon_1+\phi_n\otimes\psi_{n-1}\otimes\epsilon_2
 & (n+1)^2 & 0 \\
 & & & & \\
 \hline
 & & & & \\
 31 & SL(n)\times SL(m)\times K^*, \ n>m+1 &
 \phi_1\otimes\epsilon^a+\phi_1\otimes\psi_1\otimes\epsilon^b, \ a\ne b &
 n(m+1) & 0 \\
 & & & & \\
\hline
 & & & & \\
 32 & SL(n)\times SL(m)\times K^*, \ n>m+1>2 &
 \phi_1\otimes\epsilon^a+\phi_{n-1}\otimes\psi_{m-1}\otimes\epsilon^b,
 \ a\ne -b & n(m+1) & 0 \\
 & & & & \\
\hline
 & & & & \\
 33 & SL(n)\times SL(m)\times K^*, \ n<m &
 \phi_1\otimes\epsilon^a+\phi_1\otimes\psi_1\otimes\epsilon^b, \ a\ne 0 &
 n(m+1) & 0 \\
 & & & & \\
 & & \phi_1\otimes\epsilon^a+\phi_{n-1}\otimes\psi_{m-1}\otimes\epsilon^b,
 \ a\ne 0 & & \\
  & & & & \\
 \hline
 & & & & \\
 34 & SL(n)\times SL(2)\times SL(m), \ n>2, \ m>2 &
 \phi_1\otimes\psi_1+\psi_1\otimes\tau_1 &
 2(n+m) & 0 \\
 & & & & \\
\hline
 & & & & \\
 35 & SL(n)\times SL(2)\times Sp(2m), \ n>2, \ m\ge 1 &
 \phi_1\otimes\psi_1+\psi_1\otimes\tau_1 &
 2(n+2m) & 1 \\
 & & & & \\
\hline
 & & & & \\
 36 & Sp(2n)\times SL(2)\times Sp(2m), \ n,m\ge 1 &
 \phi_1\otimes\psi_1+\psi_1\otimes\tau_1 &
 4(m+n) & 2 \\
 & & & & \\
\hline
 & & & & \\
 37 & SL(2)\times Sp(2n)\times K^* &
 \phi_1\otimes\epsilon+\phi_1\otimes\psi_1 &
 2(2n+1) & 1 \\
 & & & & \\
\hline
\end{array}
$$

\normalsize

{\bf Comments to the Tables.} The column $"G"$ contains a
reductive group $G$.  In Table 1 the linear group $\Lambda^2
SL(n)$ is the image of $SL(n)$ under the action in the second
exterior power of the tautological representation, and $S^2 SL(n)$
is the same thing with respect to the second symmetric power.

In the column "weights" the highest weights of the $G$-module are indicated.

\noindent For the group $G_1\times G_2$ the weight
$\phi\otimes\psi$ corresponds to the tensor product of simple
$G_1$- and $G_2$-modules with highest weights $\phi$ and $\psi$
respectively. The symbol $+$ denotes a direct sum of modules. If
$G^s$ is the product of several simple groups, then their
fundamental weights are denoted successively by letters $\phi_i,\
\psi_i$ and $\tau_i$. The fundamental weight of the central torus
is denoted by $\epsilon$ (for a two-dimensional torus -- by
$\epsilon_1$ and $\epsilon_2$).

\ In the column $"\dim V"$ the dimension of the module is shown.

\ In Tables 1 and 2 the column $"{\mathcal H}"$ contains the type
of the tangent algebra ${\mathcal H}$ of the generic isotropy
subgroup $H$ for our module. Here $t_1$ is the tangent algebra of
the one-dimensional central torus in $H$, and $R_k$ is the tangent
algebra of the $k$-dimensional unipotent radical of $H$. The
information of this column is taken from Elashvili's
tables~{el1}, {el2}.

\ In the last column the codimension of a generic $G$-orbit in $V$
is shown.

\subsection{Closer look at the actions in the list}

There is a nice criterion by V. Popov, which determines when the
generic orbits of a linear action are closed, namely

\begin{Th}[cf. ~{Popov}]
Let $H$ be a stabilizer in general position for a $G$-module $V$.
Then generic orbits are closed in $V$ if $H$ is a reductive subgroup of $G$.
\end{Th}

Hence we can see that in the list of examples in the previous section,
the actions numbered 0, 2, 4, 5, 6, 7, 12, 14, 16, 17, 19, 20, 21, 23,
24, 26, 28, 29, 35, 36, 37 have closed generic orbits. These
are in fact the actions such that codimension of a generic orbit is $\geq 1$.

Now, we examine these ones to see how the generic
orbits look like topologically, in each case.

\begin{Not}
We denote by $\langle \cdot , \cdot \rangle$ the standard inner product on 
$\c^n$, i.e.
for $x = (x_1, \ldots, x_n)$ and $y = (y_1, \ldots, y_n), \langle x,y \rangle = \
sum_{i=1}^{n} 
x_iy_i$. We also denote by $\omega$ the standard symplectic form on 
on $\c^{2n}$, that is for $x = (x_1, \ldots, x_{2n})$ and $y = (y_1, \ldots, y_{2n})
, \omega(x,y) = \sum_{i=1}^{n} x_iy_i - \sum_{i=n+1}^{2n} x_iy_i$.   
\end{Not}

\begin{Def}
Let $V$ be a $G$-variety. A \textit{polynomial invariant} for the action of $G$ 
is a polynomial $P$ on $V$ such that
$$ P(x) = P(g \cdot x), \forall x \in V, \forall g \in G.$$  
\end{Def}

In the following, we find the polynomial invariant(s) of the actions 
in the previous section's list with closed generic orbits. As the codimension 
of generic orbits is 1 in almost all the cases, there, basically, exists 
only one invairant $P$ in each case and the generic orbits are given by
$\{ P(x) = c \}$ for some constant $c \in \c$.

0. $G = \{e\}$ and $V = \c$. Every point in $V$ is an orbit and $P(x) = x$.

2. $G = SL(2n, \c)$ and $V = \Lambda ^2 (\c^ {2n}) = $ vector space of
skew-symmetric matrices over $\c$. $SL(2n, \c)$ acts on $V$ by
left-right multiplication $$ g \cdot m = g^t \cdot m \cdot g.$$
Polynomial invariant is $det$.

4. and 5. $G = SL(n, \c)$ and $V = Sym^2 (\c^n) = $ vector space of 
symmetric matrices over $\c$. $SL(n, \c)$ acts on $V$ by left-right
multiplication. Again the polynomial invariant is $det$.

6. and 7. $G = SO(n, \c)$ and $V = \c^n$. $SO(n, \c)$ acts on $V$ 
via natural representation. $P(x) = \langle x,x \rangle$ is the 
polynomial invariant.

12. $G = G_2$ and $V = \c^7$. This representation is the highest 
weight representation of $G_2$ corresponding to the weight $\omega_1$
(first fundamental weight). It is the smallest representation of $G_2$ 
and is called \textit{the standard representation}. One can show that 
the action of $G_2$ on $V$ preserves a non-degenerate quadratic form (cf.
~{Fulton-Harris} p. 355).

13. $G = E_6$ and $G$ acts on $V = \c^27$. This representation is the highest 
weight representation of $E_6$ corresponding to the weight $\omega_1$
(first fundamental weight). It can be shown that the invariant polynomial 
in this case is of degree $3$ (cf. ~{}). This is the only example in 
the list that its invariant has degree $3$.

14. $G = SL(2, \c) \times \c^*$ and $V = \c^2 \otimes \c \oplus \c^2 \otimes \c = 
\c^2 \oplus \c^2$. $G$ acts by $$ (g,k) \cdot (v, w) = (k \cdot g \cdot v, 
k^{-1} \cdot g \cdot w).$$
where $v, w =in \c^2$. Invariant polynomial is $P(v, w) = det(v, w)$ (
$v$ and $w$ thought of as column vectors).

16. $G = SL(n, \c), n > 2$ and $G$ acts on $V = \c^n \oplus (\c^n)^*$ 
by $$ g \cdot (v, w)= (g \cdot v, (g^{{-1}^t}) \cdot w).$$
Invariant polynomial is  $P(v, w) = \langle v, w \rangle$.

17. and 19. $G = SL(n, \c)$ and $V = \c^{n} \oplus \Lambda^2 (\c^{n})$.
$G$ acts by $$ g\cdot (v, m) = (g \cdot v, g \cdot m).$$ 
Invariant polynomial is $P(v, m) = det(m)$.

20. $G = SO(8, \c)$ and $G$ acts on $V = \c^8 \oplus (\c^8)^8$ by 
$$ g\cdot (v, w) = (g \cdot v, (g^{-1})^t \cdot w).$$ 
There are two invariant polynomials $P(v, w) = \langle v, v \rangle$ and 
$Q(v, w) = \langle w, w \rangle$.

21. $G = Sp(2n, \c) \times \c^*$ and $G$ acts on $V = \c^{2n} \oplus \c^{2n}$ 
by $$ (g, k) \. (v, w) = (k \. g \. v, k^{-1} \. g \. w)$$. 
Invariant polynomial is $P(v, w) = \omega (v, w)$.

23. $G = SL(n, \c) \times SL(n, \c)$ and $G$ acts on $V=$ vector space of 
all $n \times n$ matrices by $$(g_1, g_2) \. m = g_1 \. m \. g_2^t .$$ 
Invariant polynomial is $P(m) = det(m)$. This is the example considered
in section ~{}, Example ~{}.

24. $G = SL(2, \c) \times Sp(2n, \c)$ and $G$ acts on $V = \c^2 \otimes \c^{2n} = $
vector space of all $2 \times 2n$ matrices  
by $$(g_1, g_2) \. m = g_1 \. m \. g_2^t .$$
If $w_1$ and $w_2$ are the rows of $m$ then the invariant polynomial is $P(m) = 
\omega(w_1, w_2)$.

26. $G = SL(4, \c) \times Sp(4, \c)$ and $G$ acts on $V = \c^4 \otimes \c^4 = $ 
vector space of all $4 \times 4$ matrices, by 
$$ g\cdot (v, m) = (g \cdot v, g \cdot m).$$ 
Invariant polynomial is $P(m) = det(m)$.

28. $G = SL(n, \c) \times SL(n, \c) \times \c^*$ and $G$ acts on 
$V = \c^n \oplus (\c^n)^{\otimes 2}$ by 
$$ (g_1, g_2, k) \. (v, m) = (k \. g_1 \. v, g_1 \. m \. {g_2}^t).$$ 
Invariant polynomial is $P(v, m) = det(m)$.

29. $G = SL(n+1, \c) \times SL(n, \c) \times \c^*$ and $G$ acts on 
$V = \c^{n+1} \oplus \c^{n+1} \otimes \c^n$ by 
$$ (g_1, g_2, k) \. (v, m) = (k^n \. g_1 \. v, g_1 \. m \. {g_2}^t).$$ 
Invariant polynomial is $P(v, m) = ??$.

35. $G = SL(n, \c) \times SL(2, \c) \times Sp(2m, \c), n > 2, m \geq 1$ and 
$G$ acts on $V = \c^n \otimes \c^2 \oplus \c^2 \otimes \c^{2m}$ by 
$$ (g_1, g_2, g_3) \. (m_1, m_2) = (g_1 \. m_1 \. {g_2}^t, g_2 \. m_2 \. {g_3}^t).$$ 
Invariant polynomial is $P(m_1, m_2) = \omega(w_1, w_2)$, where $w_1$ and $w_1$
are the rows of $m_2$.

36. $G = Sp(2n, \c) \times SL(2, \c) \times Sp(2m, \c), n, m \geq 1$ and 
$G$ acts on $V = \c^{2n} \otimes \c^2 \oplus \c^2 \otimes \c^{2m}$ by 
$$ (g_1, g_2, g_3) \. (m_1, m_2) = (g_1 \. m_1 \. {g_2}^t, g_2 \. m_2 \. {g_3}^t).$$ 
There are two invariant polynomials $P(m_1, m_2) = \omega(v_1, v_2)$ and 
$Q(m_1, m_2) = \omega (w_1, w_2)$, where $v_1$, $v_2$, $w_1$ and $w_2$ are the rows 
of $m_1$ and $m_2$ respectively.

37. $G = SL(2, \c) \times Sp(2n, \c) \times \c^*$ and 
$G$ acts on $V = \c^2 \oplus \c^2 \otimes \c^{2n}$ by 
$$ (g_1, g_2, k) \. (v, m) = (k \. g_1 \. v, g_1 \. m \. {g_2}^t).$$ 
Invariant polynomial is $P(v, m) = \omega(w_1, w_2)$, where $w_1$ and $w_1$
are the rows of $m$.

\subsection{Euler characteristic of sections and the number of
critical points of generic functionals for a hypersurface defined
by a non-degenerate quadratic form.}

Let $Q(x_1, \ldots x_n)$ be a non-degenerate quadratic form in $n$
variables over complex numbers. From linear algebra we know that
with a linear change of basis we can assume $Q(x_1, \ldots , x_n)
= x_1 ^2 + \cdots + x_n^2$. Let $X = \{ x \in \c^n | Q(x) = c \}$
for some generic constant number $c$. Obviously $deg(X) = 2$ since
a plane of complementary dimension to $X$ is a line, and
generically a line intersects $X$ in two points.

As usual let $f$ be a generic linear functional on $\c^n$ and
$(f_1, \ldots, f_n)$ the coordinates of $f$ in the dual basis. By
Lagrange's multipliers, $x = (x_1, \ldots, x_n) \in X$ is a
critical point for $f_{|X}$ iff $Q(x) = c$ and there exists scalar
$\lambda$ such that $\nabla f(x) = \lambda \cdot \nabla Q(x)$.
Notice that $\nabla f = (f_1, \ldots, f_n)$ and $\nabla Q(x) =
(2x_1, \ldots, 2x_n)$. Hence $x$ is a critical point iff $$ (x_1,
\ldots, x_n) = (f_1/2\lambda, \ldots, f_n/2\lambda)$$ where
$\lambda$ is chosen such that $ Q(f_1/2\lambda, \ldots,
f_n/2\lambda) = c$. From the latter equation one obtains $\lambda
= \pm (Q(f_1/2, \ldots, f_n/2)/c)^{1/2}$. If $f$ is generic there
are two solutions for $\lambda$ and accordingly two solutions for
$x$, that is there are two critical points for $f_{|X}$.

Next we verify the formula for the Euler characteristic of
sections of $X$, it is well-known that

\begin{Prop} Let $Q(x)$ be a non-degenerate quadratic form in 
$n$ variables over $\c$. For generic $c in \c$
the hypersurface $Q^{-1}(c)$ has the homotopy type of 
a sphere of real dimension $n-1$. 
\end{Prop}

Let  $f^{-1}(d)$ be a generic level set of $f$, that is a
hyperplane in $\c^n$. We are interested in the topology of the
intersection of $X$ and this hyperplane, i.e.
$$ X \cap f^{-1}(d) = \{ (x_1, \ldots, x_n) | Q(x_1, \ldots, x_n)
= c, f_1x_1 + \cdots + f_n x_n = d \}.$$

With out loss of generality suppose $f_n \neq 0$, then $(d -
f_1 x_1 - \cdots - f_{n-1}x_{n-1}) / f_n = x_n$. The set $X \cap
f^{-1}(d)$ is then homeomorphic to

$$ S = \{ (x_1, \ldots, x_{n-1}) | Q(x_1, \ldots, x_{n-1}, (d - f_1x_1
- \cdots - f_{n-1}x_{n-1})/f_n) = c \}.$$

The equation defining $S$ is a non-degenerate quadratic form in
$x_1, \ldots, x_{n-1}$ and hence $S$ has the homotopy type of sphere
of real dimension $n-2$.
From this, one can easily verify the formula
$$ \chi(X \cap f^{-1}(d)) = \chi(X) + (-1)^{dim(X)+1} \cdot
deg(X).$$

\subsection{Euler characteristic of sections and the number of critical points
of generic functionals for the subvariety
$\{det=\text{constant}\}$ in the space of matrices.}

There is a standard non-degenerate bilinear form on the vector
space $M(n, \c)$ of square matrices namely, $(A, B) = tr(A \cdot
B)$. Let $f$ be a functional on $M(n, \c)$ then there is a matrix
$F$ such that $f(A) = (F, A), \forall A \in M(n, \c)$.

Take a number $c$ and consider the variety $X = \{A \in M(n, \c) |
det(A) = c \}$. We want to see how many critical points  $f_{|X}$
have. As usual one uses Lagrange's multipliers: $M$ is a critical
point for iff $det(M) = c$ and there exists scalar $\lambda$ such
that $\nabla f(M) = \lambda \cdot \nabla det(M)$, where gradient
is taken with respect to the bilinear form $(\cdot ,\cdot )$. Let
us compute these gradients: since $f$ is linear and $f(A) = (F,
A)$, we have $\nabla f$ is constantly equal to $F$. As for $\nabla
det$ we know that derivative of determinant at identity is trace,
i.e. $d(det)_I(V) = tr(V)$ where $I$ is the identity matrix. So
\begin{eqnarray*}
d(det)_M(V) &=& \frac{\partial}{\partial{t}}_{|t=0} det(M + tV)
\cr &=& det(M) \cdot \frac{\partial}{\partial{t}}_{|t=0} det(I +
tM^{-1}V) \cr &=& det(M) \cdot d(det)_I(M^{-1}V) \cr &=& det(M)
\cdot tr(M^{-1}V) \cr &=& tr(det(M)M^{-1} \cdot V). \cr
\end{eqnarray*}
Hence $\nabla det(M) = det(M)M^{-1}$. So $M$ is a critical point
iff $det(M) = c$ and there exists $\lambda$ such that
$$ F = \lambda \cdot det(M) \cdot M^{-1}.$$
Taking determinant of both sides and solving for $\lambda$ we
obtain
$$ det(F) / c^{n-1}= \lambda^n.$$
and hence if $F$ is generic, there are $n$ solutions for $\lambda$
and accordingly $n$ solutions for $M$, that is $f_{|X}$ has $n$
critical points.

As determinant is a degree $n$ polynomial, degree of $X$ as a
subvariety is $n$ so we see that in this case the number of
critical points and degree are the same.

\subsection{Example of $E_6$ acting on $\c^{27}$}
As one sees in the list of actions with spherical orbits, the 
$27$ dimensional standard representation of $E_6$ is a module 
with spherical orbits. This is the only example in the list that
the invariant polynomial is of degree $3$ (cubic). Let us briefly
explain how one can construct this representation and the 
invariant cubic polynomial.
 
Let $\o$ be the $8$ dimensional algebra of Cayley numbers over $\r$.
One constructs a Jordan algebra $\j$ of $3 \times 3$ Hermitian 
matrices over $\o$. This is defined as the set of matrices of the 
form  
$$x = \left[\begin{matrix}
\alpha & a & b \\
\bar{a} & \beta & c\\
\bar{b} & \bar{c} & \gamma \\
\end{matrix}
\right] .$$
with $\alpha, \beta$ and $\gamma \in \r$ and $a, b$ and $c$ in $\o$
.The product $o$ in $\j$ is
given by
$$ xoy = \frac{1}{2}(xy + yx).$$
where the product in the right hand side are usual matrix multiplication.
This algebra is commutative but not associative and satisifies the 
identity $((xox)oy)ox = (xox)o(yox)$. One can define $tr$ and $det$
functions on $\j$ analogous to $tr$ and $det$ of matrices over fields.
They are defined by
$$ tr(x) = \alpha + \beta + \gamma.$$
$$ det(x) = \alpha \beta \gamma + tr(a(c\bar{b})) - \alpha N(c) - 
\beta N(b) - \gamma N(a).$$
where $N$ denotes norm of a octonion which is sum of squares of its 
coordinates.

Note that $tr(x)$ and $det(x) in \r$. 
It is not true in general that $det(xoy) = det(x) \. det(y)$, but one can 
prove that this is the case if $x$ and $y$ generate an associative subalgebra
of $\j$. As one can readily verify, $\j$ is a $27= 3\.8 + 3$ dimensional 
vector space over $\r$ and $\det$ is a homogenuous cubic polynomial in 
$27$ variables. 

Let us consider $\j_\c$, the complexification of $\j$. It is a 
$27$ dimensional vector space over $\c$. It is well-known that
\begin{Th}
\begin{enumerate}
\item The exceptional complex algebraic group $F_4$ is the group of
automorphisms of the vector space $\j_\c$ which preserve the scalar
product $(x,y) = tr(xoy)$ and the scalar triple product $(x,y,z) = 
tr((xoy)oz)$.
\item The exceptional complex algebraic group $E_6$ is the group of
automorphisms of the vector space $\j_\c$ which preserve $det$.
\end{enumerate} 
\end{Th}

This action of $E_6$ on $\j_\c$ is called the \textit{standard representation
of $E_6$}. As one can see in the list, it is an action with spherical 
orbits. One can show that in this representation the stabilizer of 
a generic point is isomorphic to $F_4$.

Our theorem on Euler characteristic of sections in the case of orbits of 
this representation then becomes

\begin{Th}
Let $X$ be a generic orbit in the standard representation of $E_6$,
then $\chi(X) = 0$ and
$$\chi(\text{a generic hyperplane section of $X$}) = -3.$$ 
\end{Th}
\begin{proof}
Since the stabilizer of a generic orbit is $F_4$, $X$ as a homogenuous
space is $E_6/F_4$. $rank(E_6) = 6$ and $rank(F_4) = 4$ so  
$\chi(E_6 / F_4) = 0$ (cf. Theorem {}). Since $det$ is constant 
along the orbits and $det^{-1}(c)$ is connected, we see that 
$X$ is equal to $det^{-1}(c)$ for some $c \in \c$ and, as 
$det$ is a cubic homogeneous polynomial, $deg(X) = 3$. Now from 
our main theorem on Euler characteristic of section we get
$$\chi(\text{a generic hyperplane section of $X$}) 
= 0 + (-1)^{26+1} \.3 = -3.$$   
\end{proof}

\chapter{Chern classes, Euler characteristic and the number of
critical points}

In this section we will mention some basic facts about Chern
classes and then explain formulae for Euler characteristic of
hypersurfaces and the number of critical points of functions in
terms of intersection numbers of Chern classes.

\section{Basic Facts}

Let $M$ be a real differentiable manifold and $V$ a complex vector
bundle of rank $n$ over $M$. There is a natural way of associating
a sequence of characterisitc cohomology classes to $V$, namely \textit{Chern 
classes} of $V$, such that 
$c_i(V) \in H^i(M,\z), i=0, \ldots, n$ ($c_0(V)$ is defined to be $1$) and  
the following are satisfied

1. If $V$ is a trivial bundle then $c_i(V) = 0$ for $i = 1, \ldots, n$.

2. If $f: M \rightarrow N$ is a smooth map between manifolds and
$E$ a complex vector bundle over $N$ then $$c_i(f^*E) = f^*(c_i(E)).$$
where $f^*E$ is the pull-back bundle over $M$ and $f^*$ on the 
right hand side is the induced map on the cohomology.

Total Chern
class of $V$ denoted $c(V)$ is the sum $1+c_1(V)+ \cdots + c_k(V)$
of all the Chern classes.

3. If $E$ and $F$ are complex vector bundles over $M$ then $$c(E \oplus F) = 
c(E) \cdot c(F).$$
This is known as the \textit{Whitney product formula}.

There are many different ways to construct Chern classes. For construction 
of Chern classes one can refer to ~{Griffith-Harris} p.   ,
~{Bott-Tu} p. or ~{Hirtzebruch} p.

Now, let $M$ be a complex manifold. Chern classes of the tangent
bundle $TM$ are simply called Chern classes of $M$. The last Chern
class of $M$ is called the \textit{Euler class}. It has the
property that its integral is equal to $2\pi$ times that Euler
characteristic of $M$.

Now, let $M$ be a compact complex manifold. One knows that there
is a 1-1 correspondence between divisor classes and the line
bundles on $M$. If $L$ is a line bundle and $\sigma$ a section of
$L$, then the divisor class corresponding to $L$ is given by:
$$ D = \sigma \cap M. $$

\begin{Prop}[cf. ~{Griffith-Harris}]
Let $M$ be a compact complex manifold and $L$ a line bundle. The
divisor class $D$ corresponding to $M$ is the Poincare dual to the
first Chern class $c_1(L)$.
\end{Prop}

\section{Chern classes of complete intersections of hypersurfaces}

Let $D$ be a divisor on $M$ and $L$ its corresponding line bundle.
We want to represent the Chern classes of $D$ in terms of
intersection numbers of Chern classes of $L$ and $TM$. In
particular we will get a formula for the Euler characteristic of
$D$. The following proposition tells us how one can relate tangent
bundle of $D$, tangent bundle of $M$ and line bundle $M$

\begin{Prop}[cf ~{Hirzebruch}]
Let $N$ denote the normal bundle to $D$. Then we have $N_{D} =
L_{D}$ and hence $TM_{D} = TD \oplus L_{|D}$. From the Whitney
product formula we then obtain
$$ c(TM_{D}) = c(TD) \cdot c(L_{|D}).$$
\end{Prop}

Let us denote the Poincare duality map by $P: H^*(M, \z)
\rightarrow H_*(M, \z)$. Note that $c(TM_{D})$ is simply $c(TM)
\cdot P(D)$, where the dot $\cdot$ denotes the cup product of
cohomology classes. We also know that $c(L) = 1 + c_1(L) = 1 +
P(D)$. Hence
$$ c(TM) \cdot P(D) = c(TD) \cdot (1+P(D)).$$
Solving the above equation for $c(TD)$ we obtain

\begin{Prop}
Let $M$ be a compact complex manifold of dimension $n$ and $D$
divisor. Then the Chern classes of $TD$, the tangent bundle of
$D$, can be computed in terms of cup products of Chern classes of
$M$ and $P(D)$ by
$$ c(TD) = c(TM) \cdot P(D) \cdot (1+P(D))^{-1}.$$
where we interpret $(1+P(D))^{-1}$ as the Taylor series $1 - P(D)
+ P(D)^2 - \cdots$ and terms of degree higher than $n$, the
dimension of the manifold are zero. Equating degree $n$ terms in
the above we obtain the formula for the Euler class
$$ e(TD) = c_{n-1}(TD) = \sum_{i=0}^{n-1} (-1)^{n-i+1} \cdot c_i(M) \cdot P(D)^{n-i}.$$
\end{Prop}

One can generalize the above formula to intersections of several
hypersurfaces. We state the proposition without proof. The proof
follows the same lines as the previous Proposition.

\begin{Prop}
Let $D_1, \ldots ,D_k$ be $k$ transversally intersecting divisors
on $M$ and $L_1, \ldots, L_k$ their corresponding line bundles.
Let $Y = D_1 \cap \ldots \cap D_k$. One has
$$TM_{Y} = L_1 \oplus \cdots \oplus L_k \oplus TY.$$
and hence
$$ c(Y) = c(M) \cdot \prod_{i=1}^{k} P(D_i)P((1+D_i)^{-1}).$$
\end{Prop}

For example, from the above Proposition we can get a formula for
Euler characteristic of transversal intersection of two divisors,
namely
\begin{eqnarray*}
e(D_1 \cap D_2) &=& \text{sum of the degree n terms in}(P(D_1)
\cdot P(D_1) \\ \cdot (1+P(D_1)) \cdot (1+P(D_2)) \cdot c(M) \cr
&=& \sum_{2 \leq i+j \leq n} (-1)^{i+j} D_1^i \cdot D_2^j \cdot
c_{n-i-j}(M) \cr
\end{eqnarray*}

\section{Euler charactersitic of affine sections and the
number of critical points}
Now let $L$ be an ample line bundle on $M$.

Let $M \subset \c P^N$ be a smooth projective subvariety of
dimension $n$. Let $S$ be the God-given line bundle $S$ on $\c
P^N$, namely \textit{universal subbundle}. The fibre above each $x
\in \c P^N$ is simply the line through origin representing that
point in the projective space. Restriction of $S$ to $M$ gives a
line bundle $L$ on $M$. Each projective hyperplane in $\c P^N$
defines a homology class and one can show that it is in fact the
dual of the first Chern class of $S$ and hence hyperplanes
represent the divisor class corresponding to $S$. Consequently
intersection of $M$ and a hyperplane in $\c P^N$ represents the
divisor class corresponding to $L$. Let us denote the divisor
class of hyperplane sections of $M$ by $D$. Just from the
definition of degree of a subvariety we have

$$ deg(M) = D^n.$$
where the power notation in the left hand side refers to the
intersection product of homology classes.

Considering the dual of the formula for the Euler class we get

$$ \chi(D) = (-1)^{n+1}deg(M) + \sum_{i=1}^{n-1} (-1)^{n-i+1} \cdot c_i(M)
\cdot P(D)^{n-i}.$$

Now let $D' = D \cap \c^N$ be the affine part of $D$. We can also
write down a formula for Euler characteristic of $D'$. Denote by
$H$ the projective hyperplane at infinity. We then have $D = D'
\cup (D \cap H)$. By additivity of Euler characteristic we then
obtain
$$ \chi(D') = \chi(D) - \chi(D \cap H) $$
Substituting the formula for the Euler characteristic of the
intersection of two divisors we get $$ \chi(D') = \chi(D) -
\sum_{2 \leq i+j \leq n} (-1)^{i+j} D^i \cdot H^j \cdot
c_{n-i-j}(M).$$ Note that $H$ is in the same cohomology class as
$D$ thus
$$ \chi(D') = \chi(D) - \sum_{2 \leq k
\leq n} (-1)^{k} D^k \cdot c_{n-k}(M). $$ Finally substituting the
formula for $\chi(D)$ we obtain

\begin{Prop}
Let $D$ be a divisor on a projective subvariety $M \subset \c
P^N$. Let $D' = D \cap \c^N$ be the affine part of $D$. Then
$$ \chi(D') = (-1)^{n+1}deg(M) + \sum_{i=1}^{n-1} (-1)^{n-i+1} \cdot c_i(M)
\cdot D^{n-i} - \sum_{2 \leq k \leq n} (-1)^{k} D^k \cdot
c_{n-k}(M).$$ \label{Euler-Chern
}
\end{Prop}

As A.G. Khovanskii has shown in ~{Askold}, using the above
formula for Euler characteristic of affine hyperplane sections and
the knowledge of Chern classes of projective toric varieties, one
can give another proof of Bernstein's theorem. In fact, one shows
that in {Euler-Chern} and when $M$ is a projective toric
variety, all the terms except $deg(M)$ cancels out and one is left
with
$$\chi(D') = (-1)^{n+1}deg(M).$$

On the other hand, from Theorem ~{Euler-char-formula} we have
$$ \chi(D') = \chi(M \cap \c^N) + (-1)^{n+1} \cdot \mu(f, M \cap
\c^N).$$ \label{Euler-critical}
where $\mu(f, M \cap \c^N)$ is the
number of critical points of a generic linear functional $f$ on
affine part of $M$.

Note that $\chi(M \cap \c^N) = \chi(M) - \chi(M \cap H)$ and
$\chi(M \cap H) = \chi(D)$. Hence
$$ \chi(D') = \chi(M) - \chi(D) + (-1)^{n+1} \cdot \mu(f, M \cap
\c^N)$$.

Now comparing ~{Euler-critical} and ~{Euler-Chern} we can
get a formula for the number of critical points in terms of degree
of $M$ and intersection numbers of $D$ and the Chern classes of
$M$

\begin{Prop}
Let $M \subset \c P^N$ be a smooth projective subvariety, and $f$
a generic linear functional on $\c^N$. Then $\mu(f, M \cap \c^N)$,
the number of critical points of $f$ restricted to $M \cap \c^N$,
can be obtained from
\begin{eqnarray*}
\mu(f, M \cap \c^N) &=& (-1)^{n+1} (\chi(M) - 2\chi(D) +
\chi(D^2)) \cr &=& (-1)^{n+1}(c_n(X) - 2 \cdot (\sum_{i=0}^{n-1}
(-1)^{n-i+1} \cdot c_i(M) \cdot D^{n-i}) +  \\ \sum_{2 \leq k \leq
n} (-1)^{k} D^k \cdot c_{n-k}(M) \cr
\end{eqnarray*}
\end{Prop}

\end{document}